\documentclass[12pt, reqno]{amsart}
\usepackage{amsmath}
\usepackage{amssymb}
\usepackage{epsfig}
\usepackage{graphicx}
\usepackage{enumitem}
\usepackage{color}
\usepackage{fullpage}
\definecolor{shadecolor}{gray}{0.875}
\usepackage{amscd}
\usepackage{comment}
\usepackage{adjustbox}
\usepackage{multirow}
\usepackage{hyperref}

\usepackage[dvipsnames]{xcolor}

\usepackage{tikz}
    \usetikzlibrary{shapes.arrows}
    \usetikzlibrary{decorations.pathreplacing}
\usepackage{tikz-cd}
\usetikzlibrary{decorations}

\numberwithin{equation}{section}

\input xy
\xyoption{all}

\calclayout
\allowdisplaybreaks[3]

\theoremstyle{plain}
\newtheorem{prop}{Proposition}[section]
\newtheorem{Theo}[prop]{Theorem}
\newtheorem{theo}[prop]{Theorem}

\newtheorem{lemm}[prop]{Lemma}

\theoremstyle{definition}
\newtheorem{defi}[prop]{Definition}

\newtheorem*{motques}{Motivating Question}

\newtheorem{rema}[prop]{Remark}

\newtheorem{exam}[prop]{Example}

\newcounter{enumi_saved}
\newcommand{\mbb}{\mathbb}
\newcommand{\QQ}{\mbb{Q}}

\newcommand{\ZZ}{\mbb{Z}}

\newcommand{\RR}{\mbb{R}}

\newcommand{\PP}{\mbb{P}}

\def\cM{{\mathcal M}}

\def\cO{{\mathcal O}}

\def\bP{{\mathbb P}}

\def\Eff{\overline{\mathrm{Eff}}}

\def\Nef{\mathrm{Nef}}

\newcommand{\wt}{\widetilde}

\newsavebox{\sembox}
\newlength{\semwidth}
\newlength{\boxwidth}

\newsavebox{\semrbox}
\newlength{\semrwidth}
\newlength{\boxrwidth}

\makeatother
\makeatletter

\author{Roya Beheshti}
\address{Department of Mathematics \\
Washington University in St.~Louis \\
St. Louis, MO \, \, 63130}
\email{beheshti@wustl.edu}

\author{Brian Lehmann}
\address{Department of Mathematics \\
Boston College  \\
Chestnut Hill, MA \, \, 02467}
\email{lehmannb@bc.edu}

\author{Carl Lian}
\address{Department of Mathematics \\
Tufts University \\
177 College Avenue \\
Medford, MA 02155}
\email{Carl.Lian@tufts.edu}

\author{Eric Riedl}
\address{Department of Mathematics \\
University of Notre Dame  \\
255 Hurley Hall \\
Notre Dame, IN 46556}
\email{eriedl@nd.edu}

\author{Jason Starr}
\address{Department of Mathematics\\ Stony Brook University\\ Stony Brook, NY 11794 USA}
\email{jstarr@math.stonybrook.edu}

\author{Sho Tanimoto}
\address{Graduate School of Mathematics, Nagoya University, Furocho Chikusa-ku, Nagoya, 464-8602, Japan}
\email{sho.tanimoto@math.nagoya-u.ac.jp}

\title{On the asymptotic enumerativity property for Fano manifolds}

\begin{document}

\begin{abstract}
We study the enumerativity of Gromov-Witten invariants where the domain curve is fixed in moduli and required to pass through the maximum possible number of points.  We say a Fano manifold satisfies asymptotic enumerativity if such invariants are enumerative whenever the degree of the curve is sufficiently large.  Lian and Pandharipande speculate that every Fano manifold satisfies asymptotic enumerativity.  We give the first counterexamples, as well as some new examples where asymptotic enumerativity holds. The negative examples include special hypersurfaces of low Fano index and certain projective bundles, and the new positive examples include many Fano threefolds and all smooth hypersurfaces of degree $d \leq (n+3)/3$ in $\bP^n$.
\end{abstract}

\maketitle

\section{Introduction}

The problem of ``counting'' subvarieties of a fixed projective variety has a long and illustrious history (see Section \ref{sect:history}).  One of the most successful modern theories for ``counting curves'' is Gromov-Witten theory.  However, it is often the case that the virtual counts arising in Gromov-Witten theory do not answer the corresponding geometric enumerative questions. For example, the higher-genus Gromov-Witten invariants of $\bP^3$ include many positive-dimensional contributions from the boundary of the moduli space of stable maps, and therefore often fail to count smooth, embedded curves in $\bP^3$. For other targets, the question of whether the virtual counts are enumerative can be subtle. 

In this paper we focus on a particular class of Gromov-Witten (GW) invariant, where the domain curve is fixed in moduli and required to pass through the maximum number of points. Given a smooth projective variety $X$ over $\mathbb{C}$, the \emph{fixed-domain GW invariant} $\text{GW}^{X,\text{fd}}_{g,\beta,m}(p_1,\ldots,p_m)$ is a virtual count of the number of maps $f:(\widetilde{C},\widetilde{q}_1,\ldots,\widetilde{q}_m)\to X$ where the stabilization $(C,q_1,\ldots,q_m)$ of $(\widetilde{C},\widetilde{q}_1,\ldots,\widetilde{q}_m)$ is a fixed pointed curve which is general in moduli, $f(\widetilde{C})$ has class $\beta$, and the morphism $f$ satisfies $f(\widetilde{q}_i)=p_i$ for fixed general points $p_1,\ldots,p_m\in X$.  Such invariants are particularly interesting for the following reasons:
\begin{itemize}
\item  A (non-projective) K\"ahler manifold $X$ might have very few complex analytic subvarieties, but Gromov-Witten invariants with point insertions make sense for any $X$.
\item Due to the generality of the pointed curve $(C,q_1,\ldots,q_m)$ and the points $\{p_{i}\}$, we can hope to avoid the most pathological subvarieties of the moduli space of stable maps and therefore achieve transversality of intersections.  Indeed, \cite[Speculation 12]{LP23} asks whether for every Fano variety $X$, the virtual count $\text{GW}^{X,\text{fd}}_{g,\beta,m}(p_1,\ldots,p_m)$ is enumerative when the anticanonical degree is sufficiently large compared to the genus.
\item Formulas for virtual counts $\text{GW}^{X,\text{fd}}_{g,\beta,m}(p_1,\ldots,p_m)$ are much simpler in many examples than arbitrary GW invariants, see for example \cite{BP21}.
\end{itemize}
In this paper, we solve the enumerativity question for fixed-domain GW invariants of several classes of complex Fano varieties. 
We give the first  examples of Fano varieties for which the virtual counts $\text{GW}^{X,\text{fd}}_{g,\beta,m}(p_1,\ldots,p_m)$ fail to be enumerative (even in large degree), confuting \cite[Speculation 12]{LP23}.  We also give several new examples where enumerativity holds, including all smooth hypersurfaces in $\PP^n$ of degree at most $\frac{n+3}{3}$.
In practice, deciding whether $\text{GW}^{X,\text{fd}}_{g,\beta,m}(p_1,\ldots,p_m)$ is enumerative comes down to the concrete geometric question of the existence of reducible pointed curves with stabilization isomorphic to $(C,q_1,\ldots,q_m)$ and passing through the points $p_1,\ldots,p_m$. In particular, we will not need to use any properties of the Gromov-Witten counts to verify asymptotic enumerativity.

\begin{rema}
For certain targets $X$, the problem may also be considered in terms of the space of quasimaps, and the corresponding virtual counts can be related to the stable-map invariants via wall-crossing, see \cite{CFKM14,CFK14}. Many of our arguments can be adapted equally well to the quasimap setting.
\end{rema}

\subsection{Fixed-domain GW invariants}
Early computations of fixed-domain GW invariants were on Grassmannians, see for instance \cite{Ber94,BDW96}. The Vafa-Intriligator formula, proven by Siebert-Tian \cite{ST97} in large degree and by Marian-Oprea in all degrees \cite{MO07}, determines many virtual integrals of tautological classes on Quot schemes of trivial bundles on curves, and in particular determines the fixed-domain GW invariants of Grassmannians. More generally, descendent integrals on the space of stable quotients determine all Gromov-Witten invariants of Grassmannians \cite{MOP11}.

In the formulation we give below, the invariants $\text{GW}^{X,\text{fd}}_{g,\beta,m}(p_1,\ldots,p_m)$ were recently studied systematically by Buch and Pandharipande \cite{BP21} under the name ``virtual Tevelev degrees.''  The term ``Tevelev degrees'' was introduced by Cela-Pandharipande-Schmitt \cite{CPS22} after work of Tevelev \cite{Tev20} on covers of $\mathbb{P}^1$.  More recent calculations include partial results for Fano complete intersections \cite{Cel23} and blow-ups of projective spaces \cite[\S 4]{CL23}. 

We now give a precise definition.  For a smooth projective variety $X$, we let $N_{1}(X)_{\mathbb{Z}}$ denote the abelian group of numerical classes of curves and let $\Nef_{1}(X) \subset N_{1}(X)_{\mathbb{Z}} \otimes_{\mathbb{Z}} \mathbb{R}$ denote the nef cone of curves.  Fix non-negative integers $g,m$ satisfying $3g-3+m \ge 0$; this guarantees that  the moduli stack $\overline{\mathcal M}_{g,m}$ of stable curves of genus $g$ with $m$ marked points is well-defined.  Let $\beta\in N_{1}(X)_{\mathbb{Z}}$ 
be an effective curve class and let $\overline{\mathcal{M}}_{g,m}(X,\beta)$ denote the Kontsevich moduli stack of stable maps of genus $g$ and class $\beta$ with $m$ marked points on $X$.  We let $\mathcal{M}_{g,m}(X,\beta) \subset \overline{\mathcal M}_{g,m}(X,\beta)$ denote the open substack parametrizing stable maps with smooth domain.

Consider the morphism $$\phi: \overline{\mathcal{M}}_{g,m}(X,\beta) \to \overline{\mathcal{M}}_{g,m} \times X^{\times m}$$ which combines the $m$-fold evaluation map with identification of the underlying curve.  We would expect the fibers of this map to be finite precisely when the expected dimensions match, or equivalently, when
\begin{equation*}
-K_{X} \cdot \beta + \dim(X)(1-g) = m \dim(X).
\end{equation*}

\begin{defi} \label{defn-Tevelev}
  Given $X$,  a \emph{fixed-domain triple} is a triple $(g,\delta,m)$ of nonnegative
integers such that
\begin{equation*}
\delta + \dim(X)(1-g) = m\dim(X).
\end{equation*}
If furthermore we have a numerical curve class $\beta \in \Nef_{1}(X)$ of anticanonical degree $\delta$, the \emph{fixed-domain GW invariant} $\text{GW}^{X,\text{fd}}_{g,\beta,m}(p_1,\ldots,p_m)$ is the degree of $\phi_{*}[\overline{\mathcal{M}}_{g,m}(X,\beta)]^{vir}$ as a multiple of the fundamental class $\overline{\mathcal{M}}_{g,m}\times X^{\times m}$. One can more generally replace the $p_i$ with other cohomology classes on $X$, but we will work only with point incidence conditions. We will henceforth suppress the points $q_i,p_i$ from the notation and write simply $\text{GW}^{X,\text{fd}}_{g,\beta,m}$.

We say that the virtual count $\text{GW}^{X,\text{fd}}_{g,\beta,m}$ is \emph{enumerative} if a general fiber of $\phi$ consists of exactly $\text{GW}^{X,\text{fd}}_{g,\beta,m}$ smooth points, all of which correspond to maps $f:C\to X$ with irreducible (and hence smooth) domain.  Note that this includes the possibility that the general fiber of $\phi$ is empty and $\text{GW}^{X,\text{fd}}_{g,\beta,m}=0$.
\end{defi}

The virtual count $\text{GW}^{X,\text{fd}}_{g,\beta,m}$ can fail to be enumerative if a general fiber of $\phi$ contains maps $f:\widetilde{C}\to X$ where $\widetilde{C}$ is reducible, but the stable contraction of $\wt{C}$ is a general point of $\overline{\cM}_{g,m}$ (and hence irreducible). 

\begin{exam}\label{pn_example}
Suppose that $X=\bP^n$ and let $\beta$ be $e$ times the class of a line, where
\begin{equation*}
e(n+1) + n(1-g) = mn.
\end{equation*}
If $e\ge m-1$ (or equivalently, $e\le ng$) then there exist stable maps of the following form contributing to $\text{GW}^{X,\text{fd}}_{g,\beta,m}$. The domain curve consists of a smooth genus $g$ curve $C$ attached to $m-1$ rational tails. The point $p_m$ is on $C$ and each of the rational tails contains one of the points $p_1,\ldots,p_{m-1}$. The map contracts $C$ to $p_m$ and maps the $j$th tail to a rational curve on $\bP^n$ through $p_m$ and $p_j$.  (The condition that $e \geq m-1$ guarantees that we can choose the degrees of the tails so that the total degree is $e$.)  In particular, the virtual count $\text{GW}^{X,\text{fd}}_{g,\beta,m}$ fails to be enumerative.
\end{exam}

By definition, the \emph{geometric Tevelev degree} $\text{Tev}^{X}_{g,\beta,m}$ is the number of maps $f:C\to X$ in a general fiber of $\phi$ with $C$ smooth, assuming such a fiber is finite and reduced upon restriction to the open locus $\mathcal{M}_{g,m}(X,\beta)$. Geometric Tevelev degrees, while having a more transparent definition, are often more difficult to compute than the corresponding GW invariants, and at present are only fully understood for all curve classes when $X$ is a projective space \cite{L23b}. However, if the virtual count $\text{GW}^{X,\text{fd}}_{g,\beta,m}$ is enumerative, then $\text{GW}^{X,\text{fd}}_{g,\beta,m}=\text{Tev}^{X}_{g,\beta,m}$.   According to our terminology, the converse is not necessarily true: it may happen that contributions to $\text{GW}^{X,\text{fd}}_{g,\beta,m}$ from curves of reducible domain exist, but all cancel upon integration against the virtual fundamental class, see Example \ref{product_example}.

\begin{exam}\label{product_example}
Suppose that $X=\bP^1\times\bP^1$ and let $\beta$ be the class of a $(e_1,e_2)$-curve, where
\begin{equation*}
m=e_1+e_2-g+1
\end{equation*}
and furthermore $m-1\le e_1<e_2$. Note first that there are no maps $C\to\bP^1\times\bP^1$ in class $\beta$ satisfying the conditions $f(q_i)=p_i$, because the map $C\to\bP^1$ obtained by projection to the first factor would be required to satisfy $m>2e_1-g+1$ incidence conditions, which is impossible. Furthermore, we have $\text{GW}^{X,\text{fd}}_{g,\beta,m}=0$ by the GW product formula, see \cite[\S 2.2]{BP21}.

On the other hand, adapting the construction of Example \ref{pn_example} gives stable maps of $\overline{\mathcal{M}}_{g,m}(X,\beta)$  that map to a general point of $\overline{\mathcal{M}}_{g,m}\times X^{\times m}$. Thus, in our terminology, we say that $\text{GW}^{X,\text{fd}}_{g,\beta,m}$ fails to be enumerative, despite the fact that it is equal to $\text{Tev}^{X}_{g,\beta,m}$.
\end{exam}

The main question we consider is:

\begin{motques}
Let $X$ be a Fano variety and fix a numerical curve class $\beta$ and non-negative integers $g,m$ satisfying $3g-3+m \ge 0$ and $-K_{X} \cdot \beta + \dim(X)(1-g) = m \dim(X)$.  When is the fixed-domain GW invariant $\text{GW}^{X,\text{fd}}_{g,\beta,m}$ enumerative?
\end{motques}

As we saw earlier in Example \ref{pn_example}, one can only hope that $\text{GW}^{X,\text{fd}}_{g,\beta,m}$ is enumerative when the number of points $m$ (or equivalently, the anticanonical degree $\delta$ of the curve) is large compared to the genus of the underlying curve.
Lian and Pandharipande \cite{LP23} speculate that for \emph{any} Fano $X$, the virtual count $\text{GW}^{X,\text{fd}}_{g,\beta,m}$ for a fixed-domain triple $(g,\delta,m)$ is enumerative whenever $\delta$ is sufficiently large relative to $g$.  This speculation is captured by the following definition:

\begin{defi}
We say $X$ satisfies \emph{asymptotic enumerativity} for a given genus $g$ if, for all fixed-domain triples $(g,\delta,m)$ with $m$ sufficiently large (depending only on $X$ and $g$), the virtual count $\text{GW}^{X,\text{fd}}_{g,\beta,m}$ is enumerative.
\end{defi}

\begin{rema}
The notion of asymptotic enumerativity differs from that of \emph{strong} asymptotic enumerativity \cite[Definition 5]{CL23} in that we do not require in genus 0 that $\text{GW}^{X,\text{fd}}_{0,\beta,m}$ be enumerative for \emph{all} fixed-domain triples $(0,\delta,m)$.
\end{rema}

Asymptotic enumerativity has previously been verified in the following examples:
\begin{itemize}
\item $X$ is a Grassmannian  \cite{Ber94,BDW96}
\item $X=G/P$ is a homogeneous space \cite[Theorem 10]{LP23},
\item $X$ is a hypersurface of degree $d$ in $\bP^n$, where $n-1>(d+1)(d-2)$ \cite[Theorem 11]{LP23},
\item $X=\text{Bl}_{q}\bP^n$ is the blowup of $\bP^n$ in a single point \cite[Theorem 23]{CL23}
\item $X$ is a del Pezzo surface \cite[Theorem 24]{CL23}
 \end{itemize}
 
We disprove the speculation of \cite{LP23} by giving several examples of Fano varieties for which asymptotic enumerativity fails: certain Fano threefolds, hypersurfaces with a conical linear section, and Fano projective bundles over $\mathbb{P}^n$ with sufficiently negative sections. We also give several new examples where asymptotic enumerativity holds: certain Fano threefolds and many Fano hypersurfaces.  

\subsection{Negative results}

We present two families of examples of Fano varieties $X$ for which asymptotic enumerativity fails.  The first is Fano varieties which carry a ``conical'' divisor.

\begin{prop} \label{prop:introconical}
Let $X$ be a smooth Fano variety of dimension $N$, Picard rank $1$ and Fano index $r$.  Furthermore assume that:
\begin{itemize}
\item $X$ carries an irreducible divisor $D$ and a point $p \in D$ such that every point of $D$ is contained in an anticanonical line through $p$ in $D$, and
\item There is a positive integer $t\le \frac{N}{r}-1$ such that $X$ carries a family of free rational curves of anticanonical degree $tr$.
\end{itemize}
Then $X$ fails to satisfy asymptotic enumerativity for all $g>0$. If furthermore $t<\frac{N}{r}-1$
then $X$ also fails to satisfy asymptotic enumerativity for $g=0$.
\end{prop}

An example of a variety $X$ satisfying the conditions of Proposition \ref{prop:introconical} for $g >0$ is a Fermat hypersurface $X$ of degree $d\in \left[\frac{n+3}{2},n-1\right]$ in $\mathbb{P}^{n}$ when $n\ge4$. More generally, one can apply Proposition \ref{prop:introconical} to any hypersurface of degree $d\in \left[\frac{n+3}{2},n-1\right]$ which admits a hyperplane section which is a cone over a lower dimensional hypersurface of the same degree. See Examples~\ref{exam:t-general} and \ref{exam:Fermat} for more details.

Our second family of examples is certain Fano projective bundles $X$ over projective space.  When $X$ is defined by a vector bundle which is very unstable, then certain numerical classes on $X$ cannot possibly be represented by free rational curves with balanced restricted tangent bundle.  The following result ``upgrades'' this observation to a failure of asymptotic enumerativity.

\begin{prop} \label{prop:introprojbundle}
Fix positive integers $n,r$ satisfying $n > r+1$.  Consider the vector bundle $\mathcal{E} = \bigoplus_{i=0}^{r} \mathcal{O}(a_{i})$ on $\mathbb{P}^{n-r}$ where $a_{0}=0$ and $a_{i} \geq 0$ for $i \geq 1$. Suppose that
\begin{equation*}
\sum_{i=1}^{r} a_i \leq n-r,
\end{equation*}
so that $X$ is Fano, and suppose further that
\begin{equation*}
\sum_{i=1}^{r} a_i \geq 2.
\end{equation*}
Then asymptotic enumerativity fails for the projective bundle $X = \mathbb{P}_{\mathbb{P}^{n-r}}(\mathcal{E})$ in every genus.  
\end{prop}

\subsection{Positive results}

We also identify several new situations in which asymptotic enumerativity holds. The first addresses asymptotic enumerativity for a Fano threefold $X$.  We show that asymptotic enumerativity is only obstructed by the presence of certain divisors $Y \subset X$.

\begin{theo} \label{intro:fanothreefold}
Let $X$ be a smooth Fano threefold.  Suppose that there is no divisor $Y \subset X$ that is swept out by anticanonical lines.  
Then $X$ satisfies asymptotic enumerativity for every genus $g$.
\end{theo}

\begin{exam}
Theorem \ref{intro:fanothreefold} shows that smooth cubic threefolds satisfy asymptotic enumerativity.  Smooth cubic hypersurfaces of dimension $\geq 5$ satisfy asymptotic enumerativity by \cite[Theorem 11]{LP23}.  However, the case of cubic hypersurfaces of dimension $4$ is still open.
\end{exam}

Theorem \ref{theo:fanothreefold} is a stronger version of Theorem \ref{intro:fanothreefold} which includes more information about the types of curves which can violate asymptotic enumerativity.  When a smooth Fano threefold $X$ carries a divisor $Y$ swept out by anticanonical lines, it is usually possible to combine Theorem \ref{theo:fanothreefold} with a geometric argument to determine whether one can use this divisor to violate asymptotic enumerativity. 
We do not prove a precise statement in this direction but instead give several examples; see Example \ref{exam:e5failure}, Example \ref{exam:conicalquarticthreefoldexam2}, Example \ref{exam:E3_negativeexample}, and Example \ref{exam:e3nonfailure}.

\begin{rema} \label{rema:introfujita}
The Fujita invariant of a polarized variety $(Y,L)$ compares the negativity of the canonical divisor against the positivity of $L$ (see Definition \ref{defi:fujitainv}).  If a Fano variety $X$ admits a subvariety $Y$ swept out by a family of rational curves whose dimension is larger than expected then \cite[Theorem 1.1]{LTCompos} shows that $(Y,-K_{X})$ has Fujita invariant larger than $1$.  Conjecturally the reverse implication is also true.

As demonstrated in the proof of Theorem \ref{intro:fanothreefold}, the condition on $Y$ in Theorem \ref{intro:fanothreefold} is the same as requiring that the Fujita invariant of $Y$ with respect to $-K_{X}$ is at least $2$.  More generally, we expect that asymptotic enumerativity will often fail for Fano varieties $X$ which carry a subvariety with large Fujita invariant (which conjecturally admit ``too many'' rational curves).  Unfortunately it is too much to hope for an implication in either direction.  
Example \ref{exam:e3nonfailure} gives an example of a smooth Fano threefold which satisfies asymptotic enumerativity in every genus even though it carries a divisor of Fujita invariant $2$.  Conversely, Example \ref{exam:fujitafailure} shows that one can construct curves which violate asymptotic enumerativity even when no irreducible component of the curve lies on a subvariety $Y$ with Fujita invariant $>1$.
\end{rema}

Our final result builds off of \cite{LP23} to prove new examples of asymptotic enumerativity for hypersurfaces.

\begin{Theo}\label{theo:introhypersurface}
Suppose $X$ is a smooth hypersurface of degree $d$ in $\PP^n$ such that 
\begin{itemize}
    \item[(1)] $d \leq n-2$ and $X$ is general, or 
    \item[(2)] $d \leq (n+3)/3$. 
\end{itemize}
Then, $X$ satisfies asymptotic enumerativity for every genus $g$.
\end{Theo}

In particular, we obtain transversality in the geometric calculations of \cite{L23a} when $d\le(n+2)/2$ and $X$ in addition satisfies (1) or (2). Note that the linear bound (2) is a significant improvement over the quadratic bound of \cite[Theorem 11]{LP23}. 
We expect that the condition in (2) can be weakened to ``$d \leq \frac{n+1}{2}$.''  Our proof method could conceivably extend to this larger degree range, but it would require solving a conjecture concerning the dimension of the space of non-free lines for the corresponding hypersurfaces.

\subsection{History} \label{sect:history}

For a compact K\"{a}hler manifold, what effective complex analytic
cycles can we ``count'', and are those ``counts'' invariant under
deformation?  In the special case of complex projective manifolds,
this is one of the oldest problems in algebraic geometry, dating back
at least 150 years to work of Steiner, Chasles, de Jonqui\'{e}res and,
especially, Schubert \cite{EH16}. Hilbert's
search for a rigorous foundation for Schubert's (unpublished) methods
for counting cycles eventually led to the invention of the cohomology
ring \cite{Kleiman}. The modern perspective on ``counting''
cycles of complex dimension $1$ led to the invention of the
\emph{quantum cohomology ring} whose ring product is a deformation of
the usual cup product with ``corrections'' coming from
\emph{Gromov-Witten invariants}.  In particular, these are invariant
under complex deformation, and even under symplectic deformation.
However, this does not completely solve the problem of counting
cycles, since Gromov-Witten invariants allow contributions that are
classically prohibited, e.g., Gromov-Witten invariants can be
fractional and even negative whereas true ``counts'' are nonnegative
integers.  This leads to a sharper question: when are Gromov-Witten
invariants of a K\"{a}hler manifold \emph{enumerative}?

Recall a \emph{K\"{a}hler manifold} is a connected, compact
differentiable manifold $M$ whose (real) dimension is even, say $2n$,
together with a symplectic form $\omega$ (real, closed, everywhere
nondegenerate differential $2$-form) with class
$[\omega]\in H^2(M;\RR)$, and an $\omega$-compatible complex structure
$J$, i.e., the bilinear form $\omega(\bullet, J(\bullet))$ is
everywhere symmetric and positive definite.  There are many cohomology
classes associated to $(M,\omega,J)$ that are invariant under complex
deformations and even under symplectic deformations of $\omega$, e.g.,
the Chern classes $c_i(T_{M,J})$ of the holomorphic tangent bundle for
$i=0,\dots,n$, and thus also all polynomials in these classes.  In
fact, these are all of Hodge $(p,p)$-type, so potentially in the image
of the cycle class map, so perhaps arising from cycles that we can
count.  In fact, there are many K\"{a}hler manifolds where the only
connected, closed analytic cycles are $M$ and points.
By contrast, note that the real $(1,1)$-class
$[\omega]$ is almost never invariant under symplectic deformations,
yet it is always the image of an effective $\QQ$-cycle if this
$(1,1)$-class is in $H^2(M;\QQ)$, by the Kodaira embedding theorem.
We can finesse this by restricting to K\"{a}hler manifolds
$(M,\omega,J)$ such that $[\omega]$ is a specified positive rational
multiple of $c_1(T_{M,J})$, i.e., \emph{Fano manifolds}, or such that
$[\omega]$ is a specified negative rational multiple of $c_1(T_{M,J})$,
i.e., manifolds of general type with ample canonical divisor class.

In  \cite{Gromov85} Gromov showed that often there are
complex analytic subvarieties of complex dimension $1$, and these
$J$-holomorphic curves can often be ``counted'' to produce a rational
number that is independent of symplectic deformations.  Over the next decade, with much input from
Witten and other physicists, this developed into \emph{Gromov-Witten
  theory}.
The definition within 
algebraic geometry applies to every
smooth projective variety $X$ defined over an arbitrary field $k$
(of arbitrary characteristic).
The foundational papers for
Gromov-Witten invariants in algebraic geometry are \cite{Kont95},
\cite{BF97}, \cite{BM96}, \cite{Behrend97}, \cite{LT98}, and \cite{Poma15}. 
There are 
excellent surveys of the constructions in algebraic geometry: 
\cite{FP97} and \cite{CK99}.  (For constructions in symplectic side, we refer to \cite{Ruan, RT95}.)

A ``curve class'' $\beta$ on $X$ can either be interpreted as an
element in $H^{2n-2}(X)$ for an appropriate Weil cohomology theory or,
more often, as an element in the finitely generated group
$N_1(X)_\ZZ:=\text{Hom}_{\ZZ-\text{mod}}(\text{Pic}(X)/\text{Pic}^0(X),\ZZ)$.
For every integer $g\geq 0$, for every integer $m\geq 0$, and for
every curve class $\beta$, there is a (homological) algebraic cycle class (with
$\QQ$-coefficients),
$$
GW^{X,\beta}_{g,m} \in \text{CH}_d(X^{\times m} \times \mathfrak{M}_{g,m})_\QQ,
$$
where $d=d^{X,\beta}_{g,m}$ is the \emph{virtual dimension}
$$
\langle c_1(T_X),\beta\rangle + (\text{dim}(X)-3)(1-g)+m,
$$
and where $\mathfrak{M}_{g,m}$ is the Artin stack of genus-$g$,
$m$-pointed, prestable curves.  This Artin stack is smooth of
dimension $3g-3+m$.  Thus, for (cohomological) classes $\gamma_i$ of
(complex) codimensions $d_i$ on $X$ for $i=1,\dots,m$, and for a
(cohomological) class $\lambda$ of (complex) codimension $e$ on
$\mathfrak{M}_{g,m}$, the cap product pairing of the homological cycle
class $GW^{X,\beta}_{g,m}$ against the cohomological classes gives an associated
number,
$$
GW^{X,\beta}_{g,m}(\gamma_1,\dots,\gamma_m,\lambda),
$$ 
which is zero unless the sum $d_1+\dots + d_m + e$ 
equals $d^{X,\beta}_{g,m}$.  This is the typical formulation of
Gromov-Witten invariants: as a functional on the set of $m+1$-tuples
of cohomological classes.  

The cycle class $GW^{X,\beta}_{g,m}$ is the proper
pushforward with respect to a regular $1$-morphism of Artin stacks
$$
(\Phi,\text{ev}):\overline{\mathcal{M}}_{g,m}(X,\beta) \to
\mathfrak{M}_{g,m} \times X^{\times m},
$$
of a \emph{virtual fundamental class}
$$
[\overline{\mathcal{M}}_{g,m}(X,\beta)]^{\text{vir}}\in
\text{CH}_d(\overline{\mathcal{M}}_{g,m}(X,\beta))_\QQ.
$$
Here $\overline{\mathcal{M}}_{g,m}(X,\beta)$ is an Artin stack with
finite diagonal (in characteristic $0$ it is a Deligne-Mumford stack).
It parametrizes isomorphism classes of flat families of genus-$g$,
$m$-pointed stable maps $(u: \widetilde{C} \to X,\widetilde{q}_1,\dots,\widetilde{q}_m)$ to $X$ from a
connected, proper, reduced, at-worst-nodal curve $\widetilde{C}$ of arithmetic
genus $g$ with $m$ specified (rational) points $\widetilde{q}_1,\dots,\widetilde{q}_m$
contained in the smooth locus of $\widetilde{C}$, and with a specified morphism
$u:\widetilde{C}\to X$ such that the log dualizing sheaf
$\omega_{\widetilde{C}}(\underline{\widetilde{q}}_1+\dots+\underline{\widetilde{q}}_m)$ is $u$-ample.  The
$1$-morphism $\Phi$ gives the stabilization of the $m$-pointed,
prestable curve $(\widetilde{C},\widetilde{q}_1,\dots,\widetilde{q}_m)$, and $\text{ev}$ gives the
$m$-tuple $(u(\widetilde{q}_1),\dots,u(\widetilde{q}_m))$.

In particular, \emph{primary Gromov-Witten invariants} are those with
$\lambda$ equal to $1$.  When also each class $\gamma_i$ is the
Poincar\'{e} dual of the class of a closed subscheme $Z_i$ of $X$ with
pure dimension $d_i$, then the cap product above is the pushforward of
a $0$-cycle inside the closed substack
$\text{ev}^{-1}(Z_1\times \dots \times Z_m) \subset
\overline{\mathcal{M}}_{g,m}(X,\beta)$ parameterizing stable maps with
$u(\widetilde{q}_i)$ contained in $Z_i$.  When this closed substack is a disjoint
union (possibly empty) of points, then the Gromov-Witten invariant is
\emph{enumerative}: the Gromov-Witten invariant equals the number of
stable maps mapping each marked point $p_i$ into $Z_i$.

Notice that the expected dimension is negative if $\langle
c_1(T_X),\beta \rangle$ is negative and we replace $\beta$ by a
suitably positive multiple $e\beta$ of the curve class.  For this
reason, most enumerative results are proved under the hypothesis that
$\langle c_1(T_X),\beta \rangle$ is positive for all nonzero,
effective curve classes $\beta$: this is conjecturally equivalent to the condition that $X$ is Fano.

Some of the most striking examples of enumerativity are the oldest: for example,
every smooth cubic hypersurface in $\mathbb{P}^3$ contains precisely
27 lines.  One can leverage this into proving enumerativity for genus
$0$ curves on cubic surfaces (and del Pezzo manifolds more generally)
of higher anticanonical degree, e.g., \cite{GottPand}.  Combining techniques from the Mori
program with techniques (particularly deformation theory) introduced
for Gromov-Witten theory, there are now many enumerativity theorems
for genus $0$ Gromov-Witten invariants on Fano manifolds: e.g., \cite{Thomsen98,KP01,HRS04,BK13,RY16,Bourqui16,Cas04,Testa09,CS09,LTCompos,LTJAG,LT19,LT20,BLRT22,ST22,BJ22,okamura}.

Returning to an arbitrary K\"{a}hler manifold $X=(M,J)$, 
as mentioned before, the only connected complex
analytic subvarieties might be $M$ and points.  Thus, the only primary
Gromov-Witten invariants for which we can always discuss enumerativity are the
primary Gromov-Witten invariants
where every $\gamma_i$ equals the Poincar\'{e} dual of the class of a
point.  Of course the invariant is zero for degree reasons unless
$$
m\text{dim}(X) + 3g-3 +m= d = \langle c_1(T_X),\beta\rangle +
(\text{dim}(X)-3)(1-g)+m.
$$

Recall that we restrict in this paper to Gromov-Witten invariants with point insertions
where we also fix the isomorphism type of the underlying stable curve.
Just as for all Gromov-Witten invariants, there are many cases where
the virtual count $\text{GW}^{X,\text{fd}}_{g,\beta,m}$ is not enumerative.  What is remarkable is
that there are many cases where these counts are enumerative. For example, it had essentially been understood that fixed-domain Gromov-Witten invariants of Grassmannians, as computed (in various guises) in \cite{Ber94,BDW96,ST97,MO07,MOP11}, give geometric counts of curves when the anti-canonical degree is sufficiently large. This philosophy was revisited for more general targets in \cite{LP23}.
Equally remarkable, the proofs are asymptotic in the anticanonical
degree without any careful analysis of a ``base case.''  This finally
brings us to our Motivating Question: for which Fano manifolds are the virtual counts $\text{GW}^{X,\text{fd}}_{g,\beta,m}$ asymptotically enumerative?

\

\noindent
{\bf Acknowledgements:}
Part of this project has been conducted during the SQuaRE workshop ``Geometric Manin's Conjecture in characteristic $p$'' at the American Institute of Mathematics. The authors would like to thank AIM for the generous support. We thank Felix Janda for his comment and we also thank the referee for detailed comments which improve the exposition of the paper.

Roya Beheshti was supported by NSF grant DMS-2101935. Brian Lehmann was supported by Simons Foundation grant Award Number 851129.
Carl Lian was supported by an NSF postdoctoral fellowship, grant DMS-2001976, the MATH+ incubator grant ``Tevelev degrees,'' and an AMS-Simons travel grant. Eric Riedl was supported by NSF CAREER grant DMS-1945944.  Sho Tanimoto was partially supported by JST FOREST program Grant number JPMJFR212Z, by JSPS KAKENHI Grand-in-Aid (B) 23K25764, and by JSPS Bilateral Joint Research Projects Grant number JPJSBP120219935.

\section{Background}

We recall some basic definitions regarding Fano varieties and expected dimension.

\begin{defi} \label{defn-Pic1}
  A smooth projective variety $X$ of positive dimension
  is a \emph{Fano manifold of Picard
  rank one} if every ample divisor class is $\mathbb{Q}$-numerically
equivalent to a positive multiple of the first Chern class of the tangent bundle, $c_1(T_X) =
-K_X$.  The \emph{Fano index} is the largest positive integer $r$
such that $c_1(T_X)$ equals $r$ times an integral divisor class.
\end{defi}

\begin{exam} \label{ex-Pic1}
  For every $n\geq 3$, a
  smooth codimension-$c$ intersection in $\mathbb{P}^{n+c}$ of $c$
hypersurfaces of degrees $(d_1,\dots,d_c)$ is a Fano manifold of
Picard rank one if and only if $d_1 + \dots + d_c\leq n+c$, in which
case the Fano index equals $(n+c+1)-(d_1+\dots+d_c)$.
\end{exam}

For morphisms $u$ to a smooth projective variety $Z$
from local complete intersection curves $C$ of arithmetic genus $g$
such that the morphism has a
finite automorphism group, there is a deformation-obstruction theory that
gives an isomorphism of the completion of the local ring of the moduli
space of such maps at the point corresponding to $(u:C\to Z)$ with the
quotient of a power series ring of dimension $h^0$ by an ideal
generated by $h^1$ elements, for integers $h^0$ and $h^1$ that satisfy
the identity,
$$
h^0-h^1 = \text{deg}_C(u^*T_Z) + (\dim(Z)-3)(1-g).
$$
By Krull's Hauptidealsatz, this is a lower bound for the dimension of
the moduli space near $(u:C\to Z)$, and the moduli space is locally a
complete intersection (as a Deligne-Mumford stack) if this lower bound
equals the dimension.

\begin{defi} \label{defn-vdim}
The difference $h^0-h^1$ above is the \emph{virtual dimension} or \emph{expected dimension} of
the moduli space near the point $(u:C\to Z)$.
\end{defi}

The following notion plays a crucial role in this paper:

\begin{defi}
Let $X$ be a smooth projective variety.
Let $f : \mathbb P^1 \to X$ be a rational curve that is not contracted to a point.
We say $f$ is free if the restricted tangent bundle $f^*T_X$ is nef.
We say it is very free if $f^*T_X$ is ample.

When $f$ is free, the obstruction space $H^1(\mathbb P^1, N_f)$ vanishes where $N_f$ is the normal sheaf. In particular the moduli stack $\mathcal M_{0,0}(X)$ is smooth at the point corresponding to $f$.

If $f$ is a free curve, then the deformations of $f$ yield a dominant family of rational curves on $X$.  Conversely, if we have a dominant family of rational curves then for a general point $x \in X$ every curve $f$ that meets $x$ is free.  (See \cite[II.3.5 Proposition]{Kollar}.)  By \cite[II.3.13 Lemma]{Kollar} every free curve $f: \mathbb{P}^{1} \to X$ must satisfy $\deg(f^{*}T_{X}) \geq 2$.
\end{defi}

We also need the following birational invariant in the analysis of Fano threefolds:

\begin{defi} \label{defi:fujitainv}
Let $X$ be a smooth projective variety and $L$ be a big and nef $\mathbb Q$-Cartier divisor on $X$. We define the Fujita invariant $a(X, L)$ to be
\[
a(X, L) = \inf \{ t \in \mathbb R \, | \, \text{$tL+ K_X$ is pseudo-effective.} \}.
\]
It follows from \cite{BDPP} that $a(X, L)$ is positive if and only if $X$ is uniruled.
When $L$ is nef but not big, we formally set $a(X, L) = +\infty$.

When $X$ is singular, we define the Fujita invariant $a(X, L)$ by taking a smooth resolution $\beta : \widetilde{X} \to X$:
\[
a(X, L):= a(\widetilde{X}, \beta^*L).
\]
This is well-defined due to the birational invariance of $a(X, L)$ (see \cite[Proposition 2.7]{HTT15}).
\end{defi}

\subsection{Notation}
In this paper, we use Grothendieck's notation of projective bundles, i.e., $\mathbb P(V)$ parametrizes rank $1$ quotients of $V$.

We will use the following notation when we discuss Gromov-Witten invariants.  We will denote an object parametrized by $\overline{\mathcal{M}}_{g,m}(X,\beta)$ using the notation $f: (\widetilde{C}, \widetilde{q}_{1}, \ldots, \widetilde{q}_{m})  \to X$ where $f: \widetilde{C} \to X$ is a stable map and the $\widetilde{q}_{i}$ are the marked points on $\widetilde{C}$.  We will also let $(C,q_{1},\ldots,q_{m}) \in \overline{\mathcal M}_{g,m}$ denote the stabilization of the prestable curve $(\widetilde{C}, \widetilde{q}_{1}, \ldots, \widetilde{q}_{m})$.  Note that for every irreducible component $C_{j}$ of $C$ there is a unique irreducible component $\widetilde{C}_{j} \subset \widetilde{C}$ that maps birationally to $C_{j}$ under the stabilization map.  (When $\widetilde{C}_{j}$ and $C_{j}$ are isomorphic, we will often abuse notation and write $C_{j} \subset \widetilde{C}$.)  Since every marked point $q_{i}$ is contained in a unique irreducible component $C_{j}$, we can uniquely identify corresponding points $q_{i} \in \widetilde{C}_{j} \subset \widetilde{C}$  (which may be either marked points or nodes in $\widetilde{C}$).

\section{Failure of asymptotic enumerativity}

In this section we describe two types of Fano varieties $X$ not satisfying asymptotic enumerativity, giving counterexamples to \cite[Speculation 12]{LP23}.

\subsection{Conical Fano manifolds} \label{sec-conical}

For this subsection, $X$ will denote a Fano manifold of Picard rank $1$ and Fano index $r$ and we say that a rational curve $C$ on $X$ is a line if it has anticanonical degree $r$ (or equivalently, degree $1$ against an ample generator of the Picard group).

\begin{defi} \label{defn-conical}
Let $X$ be a Fano manifold with Picard rank $1$ and Fano index $r$.  For any positive integer $t$, we say that $X$ is
\emph{$t$-general} if there exists a family of free rational curves of anticanonical degree $tr$ on $X$.

We say that $X$ is \emph{conical} if also there exists an irreducible divisor $D$ in $X$
containing a point $p$, the \emph{vertex}, such that every point of $D$ is contained in a line $\ell\subset D$  through $p$.
\end{defi}

\begin{exam}
\label{exam:t-general}
Let $X\subset \mathbb{P}^{n}$ be a smooth hypersurface of degree $d\le n$ with $n\ge4$. Then, $X$ has Picard rank 1 and Fano index $n-d+1$. If $d\le n-1$ then $X$ is 1-general \cite[Proposition 2.13]{debarre}, that is, covered by free lines.  If $d=n$ then $X$ is 2-general, that is, covered by free conics \cite{Lewis85}, see also \cite[Exercise 3.8]{debarre}. More generally, a smooth Fano complete intersection $X$ of type $(d_1,\dots,d_c)$ in $\mathbb{P}^{n}$ is 1-general if the Fano index $n-(d_1+\cdots+d_c)+1$ is at least 2, and 2-general if $n=d_1+\dots+d_c$  (see \cite[Theorem 5.2] {CR19}).

A smooth Fano complete intersection $X$ as above is also conical if there exists a hyperplane section $D$ of $X$ that is itself a cone over a complete intersection of type $(d_1,\dots,d_c)$ in $\mathbb{P}^{n-2}$. For example, a Fermat hypersurface $X$ satisfies this property.
Indeed, if a Fermat hypersurface is defined by the equation
\[
x_0^d + \cdots + x_n^d = 0
\]
in $\mathbb P^n$, then the hyperplane section $D$ given by $x_0 -(-1)^{\frac{1}{d}}x_1 =0$ is a cone with the vertex $p = ((-1)^{\frac{1}{d}}:1:0:\cdots:0)$.
\end{exam}

\begin{prop} \label{prop:pertinentnotenumerative}
Let $X$ be a conical Fano manifold of dimension $N$, Picard rank $1$ and Fano index $r$.  Suppose that $X$ is $t$-general for some positive integer $t$ such that $$t\le \frac{N}{r}-1.$$
Then $X$ fails to satisfy asymptotic enumerativity for all $g>0$. If furthermore $$t<\frac{N}{r}-1,$$
then $X$ also fails to satisfy asymptotic enumerativity for $g=0$.
\end{prop}

Note that the hypothesis on $t$ can only hold if $r\le \frac{N}{2}$.

\begin{proof}
Let $(C,q_1,\ldots,q_m)\in\mathcal{M}_{g,m}$ be a general curve and let $p_1,\ldots,p_m\in X$ be general points. Let $e,m$ be any integers for which $$-K_X\cdot\beta=er=N(m+g-1).$$ For all $m$ if $g > 0$ and $m$ sufficiently large if $ g=0$, we have $$e=\frac{N}{r}\cdot m+\frac{N}{r}(g-1)\ge (t+1)m.$$

Let $(\widetilde{C},\widetilde{q}_1,\ldots,\widetilde{q}_m)$ be the nodal curve depicted in Figure \ref{fig:conical_comb}, obtained by attaching to $C$ at each $q_i$ a chain of two rational curves $R_i\cup S_i$, where $R_i$ intersects $C$ at $q_i$ and $S_i$ at $z_i$, and $\widetilde{q}_i\in S_i$ is a smooth point.

\begin{figure}[h!]
\begin{center}
\begin{tikzpicture}[xscale=0.5,yscale=0.60] [xscale=0.36,yscale=0.36]
			
       
       		\node at (-7,0) {$C$};
       	\draw [thick, black] (-6,0) to (7,0);
	\draw [thick, black] (-4,-1) to (-6,3);
		\node at (-4,-1.5) {$R_1$};
		\draw [thick, black] (-6,1) to (-4,5);
				\node at (-4,5.5) {$S_1$};
			\node at (-5,3) {$\bullet$};
			\node at (-4.25,3) {$\widetilde{q}_1$};
\node at (-4,0.5) {$q_1$};
		\node at (-6.25,2) {$z_1$};
	\node at (0,2) {$\cdots$};

\draw [thick, black] (6,-1) to (4,3);
		\node at (6,-1.5) {$R_m$};
		\draw [thick, black] (4,1) to (6,5);
				\node at (5,5.5) {$S_m$};
			\node at (5,3) {$\bullet$};
			\node at (5.75,3) {$\widetilde{q}_m$};
\node at (6,0.5) {$q_m$};
		\node at (3.75,2) {$z_m$};

\end{tikzpicture}
\caption{The curve $\wt{C}=C\cup(R_1\cup S_1)\cdots\cup (R_{m}\cup S_m)$.}\label{fig:conical_comb}
\end{center}
\end{figure}
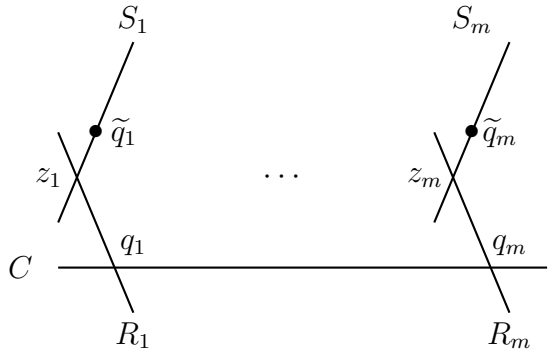

We define a stable map $f:\widetilde{C}\to X$ as follows:
\begin{itemize}
\item $f|_{C}:C\to X$ is a constant map with image $p\in D$, 
\item $f|_{S_i}:S_i\to X$ is a free curve of degree $tr$ with $f(\widetilde{q}_i)=p_i$ and $f(z_i)\in D$: such a map $f|_{S_i}$ exists because $X$ is $t$-general and of Picard rank 1, hence $D$ is ample.
\item $f|_{R_i}:R_i\to X$ is an isomorphism onto a line in $D$ with $f(q_i)=p$ and $f|_{R_i}(z_i)=f|_{S_i}(z_i)$.
\end{itemize}
As defined above, $f$ has degree $m(t+1)r\le er$, but can be modified to be a stable map of degree exactly $er$ by replacing any $f|_{R_i}$ with a multiple cover of the appropriate degree. 

Since $f$ is a point in the boundary of $\overline{\mathcal{M}}_{g,m}(X,\beta)$ over a general point of $\overline{\mathcal{M}}_{g,m}\times X^{\times m}$, $\text{GW}^{X,\text{fd}}_{g,\beta,m}$ fails to be enumerative.
\end{proof}

\begin{exam}
\label{exam:Fermat}
In particular, when $X$ is a Fermat hypersurface $X$ of degree $d$ in $\mathbb{P}^{n}$ such that $n-1 \geq d > \frac{n+3}{2}$, then asymptotic enumerativity fails in every genus.
\end{exam}

\begin{exam} \label{exam:conicalquarticthreefoldexam}
Suppose $X$ is a smooth quartic threefold with a conical hyperplane section $D$.  Proposition \ref{prop:pertinentnotenumerative} shows that $X$ does not satisfy asymptotic enumerativity for any $g>0$.
\end{exam}

\begin{exam} \label{exam:fujitafailure}
Recall that in Remark \ref{rema:introfujita} we discussed the relationship between Fujita invariants and the failure of asymptotic enumerativity.  Proposition \ref{prop:pertinentnotenumerative} shows one of the claims in Remark \ref{rema:introfujita}: the failure of asymptotic enumerativity cannot always be explained by the Fujita invariant.  More precisely, suppose that $X$ is a smooth hypersurface of degree $d$ in $\mathbb{P}^{n}$ which admits a conical hyperplane section $D$.  Furthermore suppose the degree $d$ satisfies
\begin{equation*}
\frac{n+3}{2} < d < n - 2.
\end{equation*}
The proof of Proposition \ref{prop:pertinentnotenumerative} shows that we can violate asymptotic enumerativity using curves which connect through the conical hyperplane section $D$.  However we show that $a(D,-K_{X}) < a(X,-K_{X})$.  (This should be contrasted with Theorem \ref{theo:fanothreefold} which guarantees that no such example exists amongst Fano threefolds.)

Note that $D$ is a cone over a hypersurface $Z$ of dimension $n-3$ and degree $d$ in $\mathbb{P}^{n-2}$.  The variety $\widetilde{D} = \mathbb{P}_{Z}(\mathcal{O} \oplus \mathcal{O}(1))$ is a resolution of singularities of $D$.  We let $\xi$ denote a divisor representing $\mathcal{O}_{\widetilde{D}/Z}(1)$ and $L$ denote the pullback of a hyperplane on $Z$.  Then
\begin{equation*}
K_{\widetilde{D}} = -2 \xi + (-n+d+2)L.
\end{equation*} 
The pseudo-effective cone of divisors on $K_{\widetilde{D}}$ is generated by $L$ and by $\xi - L$.  Consider divisors of the form $K_{\widetilde{D}} + a\xi$.  The assumption $d < n-2$ implies that the smallest value of $a$ for which this divisor lies in the pseudo-effective cone of $\widetilde{D}$ is $a = n-d$.  Since $-\phi^{*}K_{X} = (n+1-d)\xi$, we conclude that
\begin{equation*}
a(\widetilde{D},-\phi^{*}K_{X}) = \frac{n-d}{n+1-d} < 1 = a(X,-K_{X}).
\end{equation*}
\end{exam}

\subsection{Projective bundles}\label{subsec:projbundle}

Choose positive integers $n,r$ with $n>r+1$ and choose an $r$-tuple $(a_{1},\ldots,a_{r})$ of non-negative integers arranged in non-decreasing order.  Set $\mathcal{E} = \mathcal{O}_{\mathbb{P}^{n-r}} \oplus \mathcal{O}_{\mathbb{P}^{n-r}}(a_{1}) \oplus \ldots \mathcal{O}_{\mathbb{P}^{n-r}}(a_{r})$.  (Note that the first summand of $\mathcal{E}$ is trivial.)

Let $X = \mathbb{P}_{\mathbb{P}^{n-r}}(\mathcal{E})$ equipped with the projective bundle morphism $\pi: X \to \mathbb{P}^{n-r}$.  We let $Z$ denote the ``most rigid'' section of $\pi$ corresponding to the surjection $\mathcal{E} \to \mathcal{O}$ onto the first factor.  If we let $H$ denote the $\pi$-pullback of the hyperplane class from $\mathbb{P}^{n-r}$ and $\xi$ denote the class of the relative $\mathcal{O}(1)$, then
\begin{equation*}
\Nef^{1}(X) = \langle \xi, H \rangle \qquad \qquad \Eff^{1}(X) = \langle \xi - a_{r}H, H \rangle
\end{equation*}
Dually, if we let $\ell$ denote the class of a line in a fiber of $\pi$ and $C$ denote the class of a line in $Z$ then
\begin{equation*}
\Eff_{1}(X) = \langle \ell, C \rangle \qquad \qquad \Nef_{1}(X) = \langle \ell, a_{r}\ell + C \rangle
\end{equation*}
The anticanonical divisor is $-K_{X} = (r+1)\xi + (n-r+1 - \sum_{i=1}^{r} a_{i})H$.  In particular $X$ will be Fano if and only if $\sum_{i=1}^{r} a_{i} \leq n-r$.  

In this section we show that fixed-domain GW invariants $\text{GW}^{X,\text{fd}}_{g,\beta,m}$ fail to be enumerative for many Fano varieties $X$ given by the construction above. 

\begin{prop} \label{prop:fanoprojbundle}
Let $X = \bP_{\mathbb{P}^{n-r}}(\mathcal{E})$ be a Fano projective bundle with $n > r+1$ and with notation as above. Suppose that 
\begin{equation*}
\sum_{i=1}^{r} a_i \geq 2.
\end{equation*}
Let $g\ge0$ be any genus.

Then, there exist fixed-domain triples $(g,\delta,m)$ with $m$ arbitrarily large for which there exists a positive-dimensional family of genus $g$ reducible curves of degree $\delta$ with fixed moduli through $m$ general points. In particular, asymptotic enumerativity fails for $X$ in every genus.
\end{prop}

The construction is simple: through $m$ general points we can find a comb whose handle is a genus $g$ curve which lies in $Z$ and whose teeth are rational curves in the fibers of $\pi$.  Such combs will have higher than expected dimension.

\begin{proof}
For any $g,m\ge0$ with $3g-3+m\ge0$, define $\delta = n(m+g-1)$ so that $(g,\delta,m)$ is a fixed-domain triple. Fix general points $p_{1},\ldots,p_{m} \in X$ and let $(C,q_1,\ldots,q_m)\in \cM_{g,m}$ be a general curve. Let $(\wt{C},\widetilde{q}_1,\ldots,\widetilde{q}_m)$ be the nodal curve obtained by attaching to $C$ at $q_i$ a single rational tail $T_i$ with smooth marked point $\widetilde{q}_i$, see Figure \ref{fig:projbundle_comb}.
\begin{figure}[h!]
\begin{center}
\begin{tikzpicture}[xscale=0.5,yscale=0.60] [xscale=0.36,yscale=0.36]
			
       
       		\node at (-7,0) {$C$};
       	\draw [thick, black] (-6,0) to (6,0);
	\draw [thick, black] (-5,-1) to (-5,3);
		\node at (-5,4) {$T_1$};
\node at (-4.5,0.5) {$q_1$};
		\node at (-6,2) {$\widetilde{q}_1$};
	\node at (0,2) {$\cdots$};
	\node at (-5,2) {$\bullet$};
	\draw [thick, black] (5,-1) to (5,3);
			\node at (5,4) {$T_m$};
			\node at (6,2) {$\widetilde{q}_m$};
			\node at (4.5,0.5) {$q_m$};
		\node at (5,2) {$\bullet$};

\end{tikzpicture}
\caption{The curve $\wt{C}=C\cup T_1\cup\cdots\cup T_{m}$.}\label{fig:projbundle_comb}
\end{center}
\end{figure}
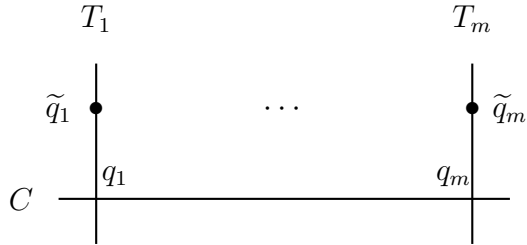

Consider the set of morphisms $f:\wt{C} \to X$ of the following form.
\begin{itemize}
\item each $T_{i}$ is mapped to a rational curve in the fiber of $\pi$ such that $f(\widetilde{q}_i)=p_{i}$ and $f(q_i)\in Z$; 
\item $C$ is mapped to $Z$, necessarily containing the points $f|_{T_i}(q_i)$
\end{itemize} 
Write $\gamma=n - r + 1 - \sum_{i=1}^{r}a_{i}$.  Let $\delta_0$ be the degree of $f|_{C}$ in $Z$ and let $m'$ denote the intersection $\xi \cdot \left( \sum_{j=1}^{m} f(T_{j}) \right)$ recording the total degree of the teeth $T_1,\ldots,T_m$ with respect to the relative hyperplane class. The condition $\delta = n(m+g-1)$ implies that these constants must satisfy
\begin{equation*}
\delta_0 \gamma + m'(r+1) = n(m+g-1).
\end{equation*}
Fix a sufficiently large positive integer $k$ and set $m=k\gamma-g+1$.  We may choose our curve $\widetilde{C}$ so that $m' = (k+1) \gamma \ge m$. We find that
\begin{equation*}
\delta_0=nk-(k+1)(r+1)=k(n-r-1)-(r+1).
\end{equation*}
Since we are assuming that $n>r+1$, the degree $\delta_{0}$ is positive for $k$ sufficiently large.

Next, consider the space parametrizing morphisms $f_0:C\to Z\cong\PP^{n-r}$ of degree $\delta_0$ for which $f_0(q_i)=f|_{T_i}(q_i)$. This space is non-empty of dimension
\begin{equation*}
\delta_{0}(n-r+1) + (n-r) - (m+g)(n-r)
\end{equation*}
so long as this number is non-negative, by \cite[Corollary 1.3]{Lar16} and the Brill-Noether theorem. Substituting our earlier value for $\delta_{0}$, this number is equal to
\begin{align*}
(k(n-r-1)-(r+1))&(n-r+1) + (n-r) - (m+g)(n-r) \\
=k[&(n-r-1)(n-r+1)-\gamma(n-r)]+B,
\end{align*}
where $B$ only depends on $n,r,g$. Thus, the dimension of this space of maps is positive as long as $k$ is sufficiently large and
\begin{equation*}
n-r+1-\sum_{i=1}^r a_i = \gamma<\frac{(n-r-1)(n-r+1)}{n-r},
\end{equation*}
which happens whenever $\sum_{i=1}^r a_i \ge 2$.

Thus, there is a positive-dimensional family of such $f:\wt{C}\to X$ with the needed properties, completing the proof.
\end{proof}

\begin{rema}
The bound $\sum_{i=1}^{r}a_i\ge2$ is sharp. When $r=1$ and $a_1=1$, $X=\PP_{\PP^{n-1}}(\cO\oplus\cO(1))$ is isomorphic to the blowup of $\PP^n$ at one point, for which asymptotic enumerativity holds by \cite[Theorem 23]{CL23}.
\end{rema}

\begin{exam}
The unique example of a Fano projective bundle of dimension $3$ that satisfies the conditions of Proposition \ref{prop:fanoprojbundle} is $X=\PP_{\PP^{2}}(\cO\oplus\cO(2))$.
\end{exam}

\section{Fano threefolds} \label{sect:fanothreefold}

In this section we study enumerativity of fixed-domain GW invariants for Fano threefolds.  Our intuition is shaped by \cite[Theorem 1.2]{BLRT22} which shows that if $X$ carries a non-dominant family of rational curves then these curves must sweep out a surface $Y \subset X$ satisfying one of the following conditions:
\begin{enumerate}
\item $Y$ is swept out by anticanonical lines, or
\item $Y$ is an exceptional divisor for a birational contraction on $X$.
\end{enumerate}
In other words, divisors $Y$ as above account for all ``unexpected'' families of rational curves.  We prove that the failure of asymptotic enumerativity for a Fano threefold is also explained by the presence of such divisors $Y$.

We will frequently use Mori's classification of the exceptional divisors $E$ on smooth Fano threefolds into five types (\cite{Mori82}):
\begin{itemize}
\item E1: an exceptional divisor for a blow-up along a smooth curve
\item E2: an exceptional divisor for a blow-up at a smooth point
\item E3: the polarized surface $(E, -K_X|_E)$ is isomorphic to $(Q, \mathcal O(1, 1))$ where $Q$ is a smooth quadric surface
\item E4: the polarized surface $(E, -K_X|_E)$ is isomorphic to $(Q, \mathcal O(1))$ where $Q$ is a quadric cone in $\mathbb{P}^{3}$
\item E5: the polarized surface $(E, -K_X|_E)$ is isomorphic to $(\mathbb P^2, \mathcal O(1))$.
\end{itemize}

In this section, we will say that a rational curve $C \subset X$ is a line (respectively conic, cubic) if $C$ has anticanonical degree $1$ (resp.~$2$, $3$).  Note that this differs from the conventions for complete intersections in $\mathbb{P}^{n}$ used in Section \ref{sec-conical} and Section \ref{sect:hypersurfaces}.  We first need the following result of \cite{LP23} which applies to stable maps with irreducible domain:

\begin{lemm}[\cite{LP23} Proposition 13 and Proposition 14] \label{lemm:curvethroughmanygenpoints}
Let $X$ be a smooth projective variety and fix $m \geq g+1$.  Suppose that $s: ({C},{q}_{1},\ldots,{q}_{m}) \to X$ is an element of $\mathcal{M}_{g,m}(X,\beta)$ that lies over a general point $((C,q_{1},\ldots,q_{m}),p_{1},\ldots,p_{m}) \in \mathcal{M}_{g,m} \times X^{\times m}$.  Then the local dimension of $\mathcal{M}_{g,m}(X,\beta)$ at $s$ has the expected value
\begin{equation*}
\dim_{[s]}\mathcal{M}_{g,m}(X,\beta) = -K_{X} \cdot \beta + (\dim(X)-3)(1-g) + m.
\end{equation*}
\end{lemm}

We also need the following estimate for stable maps with irreducible domain:

\begin{lemm} \label{lemm:curvethroughfewgenpoints}
Let $X$ be a smooth projective variety and fix $m \geq 1$. Suppose that $s: ({C},{q}_{1},\ldots,{q}_{m}) \to X$ is an element of $\mathcal{M}_{g,m}(X,\beta)$ that lies over a general point $(C,p_1)\in \cM_g\times X$, where $p_{1}$ denotes the image of ${q}_{1}$ under the first evaluation map $ev_{1}: \mathcal{M}_{g,m}(X,\beta)\to X$. 
Then the local dimension of $\mathcal{M}_{g,m}(X,\beta)$ at $s$ satisfies
\begin{equation*}
\dim_{[s]}\mathcal{M}_{g,m}(X,\beta) \leq -K_{X} \cdot \beta + \dim(X)-3 + 2g + m.
\end{equation*}
\end{lemm}

\begin{proof}
It follows from deformation theory that we have
\begin{equation*}
\dim_{[s]}\mathcal{M}_{g,m}(X,\beta) = \dim_{[s]}\mathcal{M}_{g,0}(X,\beta)+m \leq -K_{X} \cdot \beta + (\dim(X)-3)(1-g) + h^1(C, N_{s/X}) + m.
\end{equation*}
Since $s : {C} \to X$ deforms to cover $X$, we conclude that $N_{s/X}$ is generically globally generated. (See, e.g., \cite[Proposition 3.3]{LT20}) It follows from \cite[Lemma 2.8]{LRT23} that we have 
\[
h^1(C, N_{s/X}) \leq (\dim (X)-1)g.
\]
Thus our assertion follows.
\end{proof}

Our main statement for Fano threefolds relates the enumerativity of fixed-domain GW invariants with the existence of subvarieties $Y \subset X$ with large Fujita invariant.  Recall that for an $m$-pointed genus $g$ curve of anticanonical degree $d$ on a Fano threefold the fixed-domain triple condition is $d = 3m+3g-3$.

\begin{theo} \label{theo:fanothreefold}
Let $X$ be a smooth Fano threefold. Let $g\ge0$ be a fixed genus. Suppose that there exist arbitrarily large $m\ge0$ with the following property. There exists $\beta\in N_{1}(X)_{\mathbb{Z}}$ of anticanonical degree $d\le 3m+3g-3$ and an irreducible component $M_m \subset \overline{\cM}_{g,m}(X,\beta)$, such that $M_m$ is contained in the boundary of $\overline{\cM}_{g,m}(X,\beta)$, and such that the map
\begin{equation*}
\phi|_{M_m}:M_m\to \overline{\cM}_{g,m}\times X^{\times m}
\end{equation*}
is dominant. 

Then, for any $m$ sufficiently large, a general point $f:(\wt{C},\wt{q}_1,\ldots,\wt{q}_m)\to X$ of the fiber $M_{m, 0}\subset M_m$ over a general point $(C,q_1,\ldots,q_m)\in \cM_{g,m}$ has one of the following three properties.
\begin{enumerate}
\item[(i)] $f$ contracts the central genus $g$ component $C\subset \wt{C}$ to a point $p\in X$ (independent of the $q_i$), through which there exists a 1-parameter family of lines in $X$,
\item[(ii)] the central genus $g$ component $C\subset \wt{C}$ sweeps out a curve $Z\subset X$ (independent of the $q_i$) which is contained in an E5 divisor $Y\subset X$, or 
\item[(iii)] the central genus $g$ component $C\subset \wt{C}$ sweeps out a surface $Z \subset X$ with the property that $a(Z,-K_{X}|_{Z}) \ge 2$.
\end{enumerate}
\end{theo}

\cite[Proposition 4.1]{LTT18} shows that any surface $Z \subset X$ swept out by lines satisfies $a(Z,-K_{X}|_{Z}) \geq 2$.  Thus in case (i) the surface $Z$ swept out by the 1-parameter lines through $p$ will still have Fujita invariant at least 2. In case (ii), the E5 divisor $Y$ has Fujita invariant 3. Therefore, $X$ must contain a divisor of Fujita invariant at least 2 in all three cases.

\begin{proof}
Suppose that $f:(\wt{C},\wt{q}_1,\ldots,\wt{q}_m)\to X$ is a general stable map parametrized by $M_m$ over a general point of $\overline{\cM}_{g,m}\times X^{\times m}$. Then the domain of $f$ is the union of an irreducible genus $g$ curve $C$ with $b$ genus $0$ trees $T_{1},\ldots,T_{b}$.  After relabeling, we may suppose that if $T_i$ contains any marked point then it contains $\wt{q}_i$ and it meets $C$ at $q_i$.  We fix the moduli of the stabilized curve $(C,q_1,\ldots,q_m)$ throughout the proof, and assume that $\phi|_{M_m}$ remains dominant upon pullback over $[(C,q_1,\ldots,q_m)]\in\cM_{g,m}$ for arbitrarily large $m$.
That is, we assume that $M_{m, 0}$ dominates $X^{\times m}$.

We first make a series of reductions to constrain the topology of $f$. At the cost of decreasing $d$ (preserving the inequality $d\le 3m+3g-3$), we may assume that the tree $T_i$ contains the marked point $\widetilde{q}_i$, by deleting all trees $T_i$ that do not contain a marked point. We may furthermore assume that $T_i$ is a \emph{chain} of rational curves with $\widetilde{q}_i$ on the component furthest away from $C$, by successively deleting leaves on $T_i$ not containing $q_i$, also at the cost of decreasing $d$.

Now, suppose that $T_i$ has anticanonical degree at least 4. Then, deleting $T_i$ and the points $p_i\in X$ and $q_i\in C$ decreases $d$ by at least 4 and $m$ by 1, so still have $d\le 3m+3g-3$. On the other hand, this operation can be performed at most $\frac{d}{4}\le\frac{3}{4}m+O(1)$ times, where here and in the rest of the proof the implicit constant $O(1)$ is allowed to depend on $g,X$ but not on $m$. In particular, we still have $m\to\infty$ after this reduction. (We may also re-label the chains $T_i$ with the indices $i=1,2,\ldots,b$ after deleting some indices $i$.)

If any components are deleted in the process described above, then we have $d<3m+3g-3$. If the resulting curve $\wt{C}$ is irreducible (that is, $b=0$) and $m$ is sufficiently large, then by Lemma \ref{lemm:curvethroughmanygenpoints}, the map $\phi|_{M_m}$ cannot be dominant. Thus, we may assume that $b>0$.

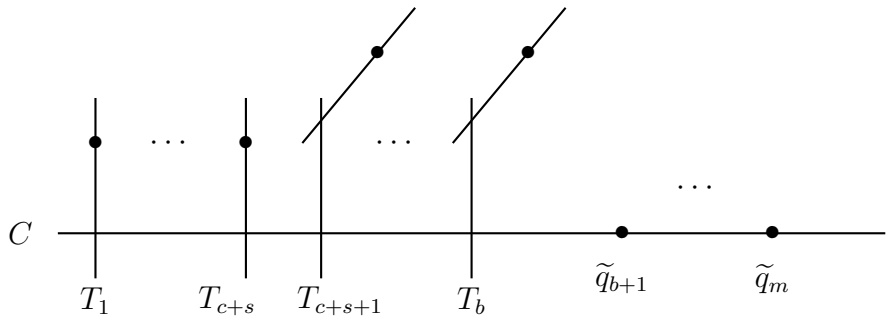
\begin{figure}[h!]
\begin{center}
\begin{tikzpicture}[xscale=0.5,yscale=0.60] [xscale=0.36,yscale=0.36]
			
       
       		\node at (-13,0) {$C$};
       	\draw [thick, black] (-12,0) to (10,0);
	\draw [thick, black] (-11,-1) to (-11,3);
	\node at (-11,-1.5) {$T_1$};
	\node at (-11,2) {$\bullet$};
	\node at (-9,2) {$\cdots$};
	
\draw [thick, black] (-7,-1) to (-7,3);
\node at (-7.5,-1.5) {$T_{c+s}$};
\node at (-7,2) {$\bullet$};

	\draw [thick, black] (-5,-1) to (-5,3);
	\node at (-4.5,-1.5) {$T_{c+s+1}$};
	\draw [thick, black] (-5.5,2) to (-2.5,5);
	\node at (-3.5,4) {$\bullet$};
	
			\node at (-3,2) {$\cdots$};

	\draw [thick, black] (-1,-1) to (-1,3);
	\node at (-1,-1.5) {$T_{b}$};
	\draw [thick, black] (-1.5,2) to (1.5,5);
	\node at (0.5,4) {$\bullet$};

	\node at (3,0) {$\bullet$};
	\node at (3,-1) {$\widetilde{q}_{b+1}$};
		\node at (5,1) {$\cdots$};
	
	\node at (7,0) {$\bullet$};
	\node at (7,-1) {$\widetilde{q}_{m}$};


\end{tikzpicture}
\caption{The curve $\wt{C}=C\cup T_1\cdots\cup T_{b}$.}\label{fig:3fold_comb}
\end{center}
\end{figure}

Because the image of each chain $T_{i}$ contains a general point of $X$, we have $\deg(T_{i}) \geq 2$ for each $i \leq b$. Indeed, the tree $T_i$ must contain a free curve of class $\beta_i$, that is, we have a dominant map $\text{ev}:M_{0,1}(X,\beta_i)\to X$. Lemma \ref{lemm:curvethroughmanygenpoints} shows that the anti-canonical degree of $\beta_i$ must be at least 2. Thus, we may assume that $\deg(T_i)\in\{2,3\}$ whenever $T_i$ is non-empty.

If $\deg(T_{i}) = 2$, then the image of $T_{i}$ is an irreducible conic through a general point.  Amongst the set of chains $T_{1},\ldots,T_{b}$, we let $T_{1},\ldots,T_{c}$ denote the trees which have anticanonical degree $2$. The stability of $f$ requires that such $T_i$ are irreducible.

If $\deg(T_{i}) = 3$, then the image of $T_{i}$ with $i \leq b$ is either the union of a conic with a line or it is an irreducible cubic. Let $T_{c+1},\ldots,T_{c+s}$ denote the chains which have anticanonical degree $3$ and contain a unique free component of degree 3. Again, the stability of $f$ requires that such $T_i$ are irreducible. Let $T_{c+s+1},\ldots,T_{b}$ denote the remaining trees which have anticanonical degree $3$, and therefore contain both a line and a free conic. Write $s'=b-c-s$. Then, by the stability of $f$ and the assumption that $T_i$ is a chain, we must have $T_i=L_i\cup C_i$, where the line $L_i$ appears between $C$ and the free conic $C_i$, which contains $\widetilde{q}_i$.

The topology of $\wt{C}$ is depicted in Figure \ref{fig:3fold_comb}. We write $d_0=d-2c-3s-3s'$ for the anticanonical degree of $C$ and $s_+=s+s'$.

As we vary the points $p_1,\ldots,p_m\in X$ in moduli, the corresponding maps $f:\wt{C}\to X$ will also deform. Let $Z\subset X$ be the subvariety of $X$ swept out by $C$. We analyze the maps that appear depending on the dimension of $Z$, and show that, in order for such stable maps $f$ through general points $p_1,\ldots,p_m$ to exist with $m$ arbitrarily large, the stable map $f$ must be of one of the three possible forms in the statement of the Theorem.

\textbf{Case 1: $\dim(Z)=0$.}

In this case, $f$ contracts $C$ to a point $p\in X$. Neither a family of free conics nor cubics through a general point of $X$ may pass through a fixed point $p\in X$, because
passing through $p$ imposes the expected number of conditions on any family of free curves.  So we have $c=s=0$ and $s'>0$. We have $T_1=L_1\cup C_1$, where $L_1\ni p$ is a line and $C_1$ is a free conic. If there is no positive-dimensional family of lines through $p$, then there are only finitely many options for $L_{1}$.  On the other hand, a general member of a family of free conics will avoid any fixed codimension $2$ locus in $X$.  Thus a free conic $C_1$ through a general point cannot meet any of the finitely many lines through $p$, a contradiction.  Thus (i) must hold.

\textbf{Case 2: $Z=X$.}

Because $C$ sweeps out all of $X$, by Lemma \ref{lemm:curvethroughfewgenpoints} the space of deformations of $f|_{C}$ that preserve the pointed moduli of $C$ has dimension at most $d_0+O(1)$.
We now estimate the dimension of the image of $M_{m, 0}$ under $\phi$. First, observe that, given a fixed choice of $f|_C$, we have:
\begin{itemize}
\item 1 parameter for the choice of a point $p_i$ on a free conic through $f(q_i)$, and
\item 2 parameters for the choice of a point $p_i$ on a free cubic through $f(q_i)$.
\end{itemize}
that is, if the restriction of $f$ to $C$ is specfied, then we get 1 (resp.~2) additional degrees of freedom for the position of $f(\wt{q}_i)=p_i$ if $T_i\ni \wt{q}_i$ is a free conic (resp.~free cubic).

Suppose instead that $T_i$ is a chain $L_i\cup C_i$. We may have a full 3-parameter family of $p_i$ arising from choices of such $T_i$. However, returning to our original set-up, note that the lines $L_{i}$ on $X$ sweep out a finite union of divisors $\{E_{j}\}$.  Since by assumption the deformations of $C$ in the family $M_{m,0}$ dominate $X$, the same is true if we consider the larger sublocus of $M_{m}$ where we fix $C$ but allow the points  $(q_{1},\ldots,q_{m})$ to vary in moduli on $C$.  In this larger family, we see that forcing $q_{i}$ to map to one of the divisors $E_{j}$ imposes one condition on the family of maps $f|_C$.  Since our original curve $(C,q_{1},\ldots,q_{m})$ is general in moduli, the conditions imposed at all such $q_i$ are independent.
Therefore, we conclude that the dimension of the image of $M_{m,0}$ under $\phi$ is at most
\begin{align*}
d_0+O(1)+c+2s+2s' &=d-(c+s_+)+O(1)\\
&=3m-b+O(1).
\end{align*}

On the other hand, in order for $M_{m, 0}$ to dominate $X^{\times m}$, we need
\begin{equation*}
3m-b+O(1) \ge 3m,
\end{equation*}
hence $b\le O(1)$. In particular, $m-b\ge g+1$ if $m$ is sufficiently large. Then, by Lemma \ref{lemm:curvethroughmanygenpoints}, it follows that $f|_{C}$ moves in a family of the expected dimension of $d_0-3g+3$ (we have subtracted $3g-3+m$ from the right hand side of the formula of  Lemma \ref{lemm:curvethroughmanygenpoints} to account for the fact that we have fixed the moduli of $C$). We now repeat the calculation of the dimension of the space of deformations of $f$, with this more precise estimate. The dimension of the image of $M_{m,0}$ under $\phi$ is at most
\begin{align*}
(d_0-3g+3)+c+2s_+&=(d-3g+3)-(c+s_+)\\
&=(d-3g+3)-b,
\end{align*}
which is $b$ less than the expected dimension.  Thus, the deformations of $f$ cannot dominate $X^{\times m}$, which has dimension $3m>(d-3g+3)-b$.

\textbf{Case 3: $\dim(Z)=1$.}

As in the case where $C$ is contracted, we need $c=0$ since a free conic through a general point does not meet a fixed curve $Z\subset X$. Also, we need $b=m$, or else $C$ sweeps out all of $X$. Thus
\begin{equation*}
d_0= d-3m\le O(1),
\end{equation*}
and in particular, the space of deformations of $f|_{C}$ has bounded dimension $O(1)$.  First suppose that a general point of $Z$ is not contained in a $1$-parameter family of lines. Then arguing in the previous case, the space of deformations of $f$ has dimension at most 
\begin{equation*}
O(1)+2s_+=O(1)+2m.
\end{equation*}
Thus this space of deformations cannot dominate $X^{\times m}$.

On the other hand, if a general (and hence, every) point of $Z$ is contained in a $1$-parameter family of lines, then these lines must sweep out a surface $Y$ carrying a $2$-dimensional family of lines.  By \cite[Lemma 4.4]{BLRT22} $Y$ is an E5 divisor, and we are now in case (ii) of the Theorem.

\textbf{Case 4: $\dim(Z)=2$.}

Let $\psi: \widetilde{Z} \to X$ denote the composition of a resolution of $Z$ with the inclusion.  Let $\widetilde{f}: C\to \widetilde{Z}$ be the map obtained by strict transform from a general point of $M_{m, 0}$ over $(C,q_1,\ldots,q_m)$. Then, by Lemma \ref{lemm:curvethroughfewgenpoints}, the map $\wt{f}$ moves in a family of dimension at most
\begin{equation*}
-K_{\widetilde{Z}}\cdot\wt{f}_{*}[C] + O(1).
\end{equation*}

Suppose for a contradiction that we have $a(\wt{Z},\psi^{*}(-K_{X})) < 2$.  Since deformations of $C$ sweep out $Z$ we have
\begin{equation*}
(K_{\widetilde{Z}} - a(\wt{Z},\psi^{*}(-K_{X}))\psi^{*}K_{X}) \cdot \wt{f}_{*}[C] \geq 0.
\end{equation*}
Choose $\epsilon > 0$ so that $2 - \epsilon > a(\wt{Z},\psi^{*}(-K_{X}))$.  Then for any $O(1)$ constants we find that in sufficiently large degrees
\begin{equation*}
-K_{\widetilde{Z}}\cdot\wt{f}_{*}[C] + O(1)\le (2-\epsilon)d_0 + O(1)
\end{equation*}
where $d_0=-K_X\cdot f_{*}[C]$ as before. The same estimate also gives an upper bound on the dimension of the space of deformations of $f|_{C}:C\to Z$.

Consider now the space of deformations of $f$. As in case 2, the free conics and cubics add 1 and 2 parameters, respectively. We claim that reducible chains of the form $T_i=L_i\cup C_i$ also add 2 parameters. First, suppose that $f(q_i)$ is contained in a 1-parameter family of lines. Then, because $f(q_i)$ is a general point of the divisor $Z$, and a family of lines cannot cover $X$, this 1-parameter family of lines must be contained in $Z$, and thus sweeps out $Z$. Therefore, the lines have anticanonical degree at least 2 on $\widetilde{Z}$, which contradicts the assumption that $a(Z,-K_{X}|_Z) < 2$. Thus, only finitely many lines pass through $f(q_i)$. Then, there is a 1-parameter family of free conics $C_i$ incident to $L_i$, and there is a 1-parameter family of points on $C_i$. In total, the chains $T_i$ also add 2 parameters to the space of deformations of $f$. 

We conclude that the space of deformations of $f$ has dimension at most
\begin{align*}
(2-\epsilon)(d-2c-3s_+)+O(1)+c+2s_+&\le (2-\epsilon)(3m-2c-3s_+)+c+2s_{+}+O(1)\\
&=(6-3\epsilon)m-(3-2\epsilon)c-(4-3\epsilon)s_{+}+O(1)\\
&=(3-\epsilon)m-(1-\epsilon)s_{+}+O(1) \\
& <3m,
\end{align*}
for $m$ sufficiently large.  In the next-to-last line, we have used the equality $c+s_{+} = b = m$.  In the last line, we have used the fact that $(1-\epsilon)s_{+}\ge0$ if $\epsilon\le 1$, whereas if $\epsilon>1$, then $(3-\epsilon)m-(1-\epsilon)s_+=2m+(1-\epsilon)(m-s_+)<2m$. 
Thus, $M_{m, 0}$ cannot dominate $X^{\times m}$.
\end{proof}

\begin{proof}[Proof of Theorem \ref{intro:fanothreefold}:]
Let $X$ be a smooth Fano threefold, and assume first that there is an infinite sequence of positive integers $m$ and nef curve classes $\beta$ on $X$ such that $(g,-K_{X} \cdot \beta,m)$ is a fixed-domain triple and a sequence of irreducible components $M_m \subset \overline{\mathcal{M}}_{g,m}(X,\beta)$ such that a general fiber of the morphism $M_m \to \overline{\mathcal{M}}_{g,m} \times X^{\times m}$ consists only of maps with reducible domains.  By Theorem \ref{theo:fanothreefold} $X$ admits a divisor $Y$ with $a(Y,-K_{X}|_{Y}) \geq 2$.

According to \cite[Theorem 4.1]{BLRT22}, if $X$ is a smooth Fano threefold and $Y \subset X$ is a surface satisfying $a(Y,-K_{X}|_{Y}) > 1$, then either $Y$ is swept out by anticanonical lines or $Y$ is an exceptional divisor on $X$.  The exceptional divisors with Fujita invariant $\geq 2$ have types E1, E3, E4, E5 and each admits a dominant family of anticanonical lines.

Therefore, if $X$ is a smooth Fano threefold without a divisor $Y$ as in the statement of Theorem \ref{intro:fanothreefold}, then a general fiber of the morphism $\overline{\mathcal{M}}_{g,m}(X,\beta) \to \overline{\mathcal{M}}_{g,m} \times X^{\times m}$ consists only of maps with irreducible (hence smooth) domains, for any fixed-domain triple $(g,-K_{X} \cdot \beta,m)$ with $m$ sufficiently large compared to $g$. By \cite[Proposition 13 and Proposition 14]{LP23}, the generic fiber of $\phi$ is furthermore reduced of dimension 0. It follows that $X$ satisfies asymptotic enumerativity.
\end{proof}

When $X$ is a smooth Fano threefold which carries a divisor swept out by lines it is usually straightforward to tell directly whether one can obstruct asymptotic enumerativity using stable maps consisting of rational curves attached to a central genus $g$ curve $C$ as described by Theorem \ref{theo:fanothreefold}. 

\begin{exam} \label{exam:e5failure}
Consider the Fano threefold $X = \mathbb{P}_{\mathbb{P}^{2}}(\mathcal{O} \oplus \mathcal{O}(2))$.  The rigid section $Z$ of the $\mathbb{P}^{1}$-bundle $X \to \mathbb{P}^{2}$ is an E5 divisor with Fujita invariant 3.  Proposition \ref{prop:fanoprojbundle} demonstrated the failure of asymptotic enumerativity on $X$ using the existence of the divisor $Z$.  As predicted by Theorem \ref{theo:fanothreefold}, the ``bad'' curves are obtained by attaching $m$ fibers of the projective bundle to a genus $g$ curve in $Z$.\end{exam}

\begin{exam}\label{exam:conicalquarticthreefoldexam2}
Let $X\subset\bP^4$ be a conical quartic threefold. Then $X$ fails to satisfy asymptotic enumerativity by Example \ref{exam:conicalquarticthreefoldexam}. By definition the conical divisor of $X$ contains a 1-parameter family of lines through the cone point $p$. As predicted by the proof of Theorem \ref{theo:fanothreefold}, asymptotic enumerativity is violated by curves $\wt{C}$ for which the central genus $g$ curve $C$ is contracted to $p$.
\end{exam}

\begin{exam}\label{exam:E3_negativeexample}
Let $X$ be the Fano threefold $\bP_{\bP^1\times\bP^1}(\mathcal \cO \oplus \mathcal \cO(1,1))$ and let $\pi:X\to\bP^1\times\bP^1$ be the projection. Let $Z\subset X$ be the section of $\pi$ corresponding to the projection $\cO \oplus \mathcal \cO(1,1)\to\cO$. Then $Z$ is an E1 divisor with Fujita invariant 2.  However, no point on $Z$ (or $X$) is contained in a 1-parameter family of lines.

Repeating the contruction of \S\ref{subsec:projbundle}, for any positive integer $m$ with $m+g\equiv1\pmod{2}$, we construct reducible curves of genus $g$ with general moduli through $m$ general points of $p_i\in X$ of anticanonical degree $3m+g-1$ as follows. Attach a $\left(\frac{m+g-1}{2},\frac{m+g-1}{2}\right)$-curve of genus $g$ in $Z\cong\bP^1\times\bP^1$ at the points $q_i\in C$ to rational tails mapping isomorphically to the fibers of $\pi$ over the $p_i$. Such a curve with general moduli in $Z$ exists, for example, by \cite[Corollary 1.3]{Lar16}.  
If $g>0$, then one may replace the rational tails with multiple covers of arbitrary degree to obtain a map of degree exactly $3m+3g-3$. In this way, $X$ fails to satisfy asymptotic enumerativity for any $g>0$.
\end{exam}

The following example shows that a Fano threefold can satisfy asymptotic enumerativity in every genus even if it carries a divisor with Fujita invariant $2$.

\begin{exam} \label{exam:e3nonfailure}
Let $X$ be the blow-up of $\mathbb{P}^{3}$ along the intersection $Z$ of a smooth quadric and a smooth cubic.  Then $X$ carries two divisors with Fujita invariant $2$: the exceptional divisor $E$ for the blow-up (of type E1) and the strict transform $Q$ of the quadric (of type E3).  However no point of $X$ is contained in a $1$-parameter family of lines.

We claim that $X$ satisfies asymptotic enumerativity in every genus.  By Theorem \ref{theo:fanothreefold}, it suffices to show that asymptotic enumerativity cannot be violated by families of reducible curves where the central genus $g$ curve $C$ deforms to sweep out $Q$ or $E$.

First we show that it is not possible to obstruct asymptotic enumerativity using a family of reducible curves $\wt{C}$ such that the genus $g$ component $C$ sweeps out $Q$.  The only dominant family of anticanonical conics $T$ on $X$ is given by the strict transforms of lines meeting $Z$ twice. For every such curve $Q \cdot T = 0$ and thus a general conic cannot meet $Q$.  To finish the argument we repeat the computation in Case 4 of the proof of Theorem \ref{theo:fanothreefold}.  In the notation of this proof, we have shown that $c = 0$ so that $s_{+} = m$.   The computation in this proof shows that curves of this type deform in dimension at most $2m + O(1)$.  Since the dimension of $X^{\times m}$ is $3m$, we see that asymptotic enumerativity cannot be violated using such curves.

Next we show that it is not possible to obstruct asymptotic enumerativity using a family of reducible curves $\wt{C}$ such that the genus $g$ component $C$ sweeps out $E$.  Since $Z$ has genus $4$, it does not receive a map from a general curve of genus $g$ for any $g$.  Thus the only way to violate enumerativity using maps to $E$ is if $C$ is a multiple cover of a fiber of $E \to Z$ or is contracted to a point.  Note that such a curve can meet at most one general anticanonical conic.  Repeating the computation in Case 4 of the proof of Theorem \ref{theo:fanothreefold}, we see that $s_{+} \geq m-1$ so that curves of this type deform in dimension at most $2m + O(1)$.  Thus asymptotic enumerativity cannot be violated using such curves.
\end{exam}

\section{Hypersurfaces} \label{sect:hypersurfaces}

The goal in this section is to prove the following theorem describing enumerativity of fixed-domain GW invariants for certain hypersurfaces.

\begin{Theo}\label{hypersurface-main}
Suppose $X$ is a smooth hypersurface of degree $d$ in $\PP^n$ such that 
\begin{itemize}
    \item[(1)] $d \leq n-2$ and $X$ is general, or
    \item[(2)] $d \leq (n+3)/3$. 
\end{itemize}
Then $X$ satisfies asymptotic enumerativity for every genus $g$.
\end{Theo}

The strategy of proof is as follows. It is enough to establish an upper bound, independent of $g$, on $h^1$ of the restricted tangent bundle for a stable map of genus $g$ to $X$ with general stabilized domain passing through $m$ points.  (See Lemma \ref{lemm:assumingh1} for a precise statement.) We obtain such a bound by degeneration. When $g=0$, we replace the general domain curve with a maximally degenerate chain of rational curves, see Lemma \ref{chain-rational} and Proposition \ref{prop:h1_bounded}, using results established in \S\ref{sect:chains} to obtain bounds on $h^1$. The case of arbitrary genus is finally reduced to the case of genus $0$ by replacing a general pointed curve of genus $g$ by one with a rational tail containing all of the marked points.

In this section, we will say that a rational curve $C \subset X$ is a line or conic if $C$ is a line or conic in the ambient projective space.

\subsection{Chains of curves on hypersurfaces}\label{sect:chains}

Let $X$ be a smooth hypersurface of degree $d$ in $\PP^n$. We assume throughout that $d\le n$ and that $n\ge4$. By the parameter space of chains of rational curves of degree $e$ with $t$ components, we mean 
the closure of the locus in $\overline{M}_{0,0}(X,e)$ parametrizing stable maps $f: C \to X$ of degree $e$ such that $C$ is a chain $C=C_1 \cup \dots\cup C_t$ of $t$ 
smooth rational curves such that $C_i$ intersects $C_{i-1}$ and $C_{i+1}$.  The {\em{expected dimension}} of the parameter space of chains of rational curves on $X$ of degree $e$ with $t$ components is
$$e(n+1-d)+n-3-t.$$
When $X$ is general and $d\leq n-2$, by \cite[Theorem 3.3]{RY16}, the space of rational curves of a given degree through any point of $X$ has the expected dimension $e(n+1-d)-2$, so every component of the parameter space of chains of rational curves in $X$ has the expected dimension.

\begin{prop}\label{prop:nonfree_bound}
If $X \subset \PP^n$ is a smooth hypersurface of degree $d$, then the space of non-free lines on $X$ has dimension at most $n+d-5$.
\end{prop} 

\begin{proof}
Suppose $S$ is a family of non-free lines sweeping out an irreducible subvariety $Y$ in $X$. Let $l$ be a general line parametrized by $S$. Then
by \cite[Theorem 2.4]{BR21}, $\dim S \leq \dim Y +n-3-h^0(N_{l/X}(-1)) \leq 2n-5-h^0(N_{l/X}(-1))$. Since $l$ is not free, the proof of \cite[Theorem 2.4]{BR21} shows that $h^0(N_{l/X}(-1)) \geq n-d$ and we get the desired result. 
\end{proof}

    \begin{prop}\label{chain_conics}
Let $X \subset \PP^n$ be any smooth hypersurface of degree $d\leq (n+3)/3$. Then every irreducible component of the space parametrizing chains of lines and free conics has the expected dimension.
\end{prop}

\begin{proof} 
We argue by induction on the number of irreducible components $t$ in the chain.  If $t=1$ then we either have a line or a free conic.  By \cite[Theorem 1.3]{BR21} the space of lines has the expected dimension $2n-d-3$ for every smooth hypersurface of degree $\leq (n+4)/2$.  Thus in both cases the dimension is the expected dimension.

For the induction step, suppose the statement holds for $t' < t$.  Write $\widetilde{C}=C_1 \cup \dots \cup C_t$ for the irreducible components of $\widetilde{C}$.  First suppose that $C_t$ is a free conic or a free line.  By induction, the deformation space of the chain of curves $C_{1} \cup \ldots \cup C_{t-1}$ has the expected dimension.  Since we are attaching a free curve to this chain, the total family also has the expected dimension. Indeed, the deformation space of the free curve has the expected dimension, and the requirement that the free curve be incident to the chain at any fixed point also imposes the expected number of conditions on the free curve.

Next suppose $C_t$ is a non-free line. There are two cases:
\begin{itemize}
\item $C_{t-1}$ is a line: In this case the space of deformations of $C_1 \cup \dots \cup C_{t-2}$ has dimension at most $(e-2)(n+1-d)+n-3-(t-2)$  by our induction hypothesis, and the space of non-free lines has dimension $\leq n+d-5$ by Proposition \ref{prop:nonfree_bound}.  Since we are assuming that $d \leq (n+3)/3$, the total dimension is at most
\begin{align*}
&[(e-2)(n+1-d)+n-3-(t-2)]+1+[n+d-5]+1 \\
\leq \quad&e(n+1-d)+n-3-t,
\end{align*}
where the terms in brackets come from the deformations of the chain and the component $C_{t-1}$, respectively, and the additional summands $+1$ from the possible degrees of freedom from the points to which these components are attached.

\item $C_{t-1}$ is a free conic: The argument is similar to the previous case.  First note that the space of free conics  through any two fixed points of $X$ has dimension $\leq n-d+1$.  Indeed, $N_{C/X}$ is a subbundle of $N_{C/\bP^n}=\cO(4) \oplus \cO(2)^{n-2}$ of degree $2n-2d$ with no negative summand.  Thus $N_{C/X}(-2)$ has either: 1) at most 
$n-d$ summands of degree $0$ and all the other summands negative or 2) at most $n-d-2$ summands of degree $0$, $1$ summand of degree $2$, and all the rest negative. In both cases $H^{0}(C,N_{C/X}(-2))$ is at most $n-d+1$.

By our induction hypothesis the space of deformations of $C_1 \cup \dots \cup C_{t-2}$ 
has dimension at most $(e-3)(n+1-d)+n-3-(t-2)$.  By choosing the connection points of the free conic $C_{t-1}$ with the curves $C_{t-2}, C_{t}$, we get that the dimension of the space of such chains is at most
\begin{align*}
&[(e-3)(n+1-d)+n-3-(t-2)]+1+[n+d-5]+1+[n-d+1] \\
\leq\quad & e(n+1-d)+n-3-t.
\end{align*}
\end{itemize} 

\end{proof}

\subsection{Asymptotic enumerativity for hypersurfaces}

We build up the proof of Theorem \ref{hypersurface-main} in several steps.  Consider the map 
$$\phi: \overline{\cM}_{g,m}(X,e) \to \overline{\cM}_{g,m}\times X^{\times m}.$$

\begin{lemm} \label{lemm:assumingh1}
Let $X$ be a smooth hypersurface of degree $d$ in $\PP^{n}$ and fix a genus $g$.  

Suppose that there is a constant $A$ depending only on $X$ and $g$ such that, for any 
positive integers $m,e$ with $e(n+1-d)=(n-1)(m+g-1)$, for a general $([C], p_1, \dots, p_m) \in\overline{\cM}_{g,m}\times X^{\times m}$, and for any $(\wt{C}, f, \widetilde{q}_1, \dots, \widetilde{q}_m)$ 
in the fiber of $\phi$ over $([C],p_1,\dots,p_m)$, we have $h^1(\wt{C},f^*T_X)\le A$. 

Suppose further that $X$ satisfies one of the hypotheses (1), (2) of Theorem \ref{hypersurface-main}. Then, Theorem  \ref{hypersurface-main} holds for $X$ and genus $g$. 
\end{lemm}

\begin{proof}
Suppose otherwise. Then, for arbitrarily large $m$ and $e$ and for general $([C], p_1, \dots, p_m) \in\overline{\cM}_{g,m}\times X^{\times m}$, there exist $(\wt{C}, f, \widetilde{q}_1, \dots, \widetilde{q}_m)$
in the fiber of $\phi$ such that $\wt{C}$ is reducible. Indeed, when $\wt{C}$ is irreducible, our assertion follows from \cite[Proposition 13 and Proposition 14]{LP23}.
It follows that $\widetilde{C}$ has a component isomorphic to $C$, to which trees of rational curves, each containing at most one marked point, are attached. Assume without loss of generality that the first $t$ marked points are on the attached trees and the last $m-t$ marked points are on $C$, see Figure \ref{fig:hypersurfacecomb}.
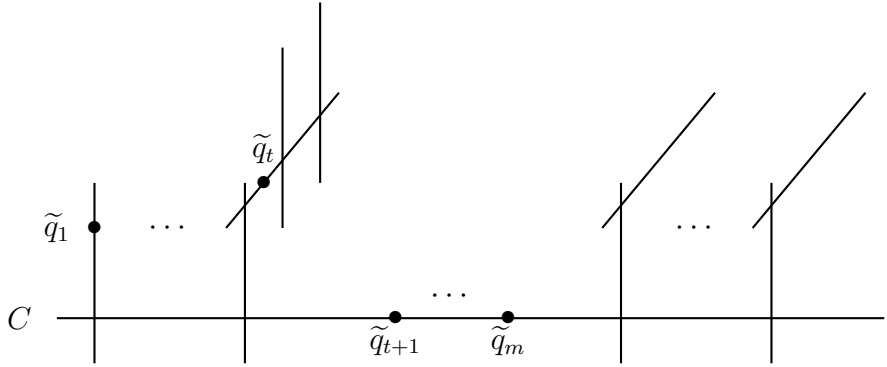
\begin{figure}[h!]
\begin{center}
\begin{tikzpicture}[xscale=0.5,yscale=0.60] [xscale=0.36,yscale=0.36]
			
       
       		\node at (-13,0) {$C$};
       	\draw [thick, black] (-12,0) to (10,0);
	\draw [thick, black] (-11,-1) to (-11,3);

	\node at (-11,2) {$\bullet$};
	\node at (-12,2) {$\widetilde{q}_1$};
	\node at (-9,2) {$\cdots$};
	
\draw [thick, black] (-7,-1) to (-7,3);
\draw [thick, black] (-7.5,2) to (-4.5,5);
\node at (-6.5,3) {$\bullet$};
\node at (-6.5,3.75) {$\widetilde{q}_t$};
\draw [thick, black] (-6,2) to (-6,6);
\draw [thick, black] (-5,3) to (-5,7);

\node at (-3,0) {$\bullet$};
\node at (-3,-0.5) {$\widetilde{q}_{t+1}$};
\node at (0,0) {$\bullet$};
\node at (0,-0.5) {$\widetilde{q}_{m}$};
\node at (-1.5,0.5) {$\cdots$};

	\draw [thick, black] (3,-1) to (3,3);
	\draw [thick, black] (2.5,2) to (5.5,5);
		\node at (5,2) {$\cdots$};
	
		\draw [thick, black] (7,-1) to (7,3);
		\draw [thick, black] (6.5,2) to (9.5,5);

\end{tikzpicture}
\caption{The curve $\wt{C}$, obtained from $C$ by attaching rational trees.}\label{fig:hypersurfacecomb}
\end{center}
\end{figure}

Let $\iota:C\to \wt{C}$ be the natural closed embedding. From the surjection $f^{*}T_X\to \iota_{*}\iota^{*}f^{*}T_X$ of sheaves on $\wt{C}$, we have 
\begin{equation*}
h^1(C,f^{*}T_X|_{C})=h^1(\wt{C},\iota_{*}\iota^{*}f^{*}T_X)\le h^1(\wt{C},f^*T_X) \le A.
\end{equation*}

First, suppose that $d<(n+3)/2$. Then, as $X$ has Fano index at least $n/2$, the hypothesis (ii) from \cite[Proposition 22]{LP23} is satisfied. Furthermore, because $h^1(C,f^{*}T_X|_{C}) \le A$, we have condition $(\star)_g$, hence condition $(\star\star)_g$, see \cite[Definition 18]{LP23}, at least upon restriction to $f$ over a general point of $\phi$. Therefore, \cite[Proposition 22]{LP23} (which only requires working over a general point of $\phi$) applies, and the conclusion of Theorem \ref{hypersurface-main} follows. In particular, if $X$ has degree $d \leq (n+3)/3$, then we obtain the conclusion. 

Now, suppose that $d\le n-2$ and that $X$ is general. Let $C'$ be the connected curve obtained by successively deleting components on the attached trees which do not contain a marked point and whose removal does not make the tree disconnected. Let $e_i$ be the degree of $f$ restricted to the irreducible component containing the $i$th marked point for $1 \leq i \leq t$. Let $C''$ be the curve obtained by removing the components with marked points on the attached trees in $C'$ and suppose the total degree of $f|_{C^{\prime\prime}}$ is $e^{\prime\prime}$. The curves $C',C''$ are depicted in Figure \ref{fig:hypersurfacecomb_deletion}.
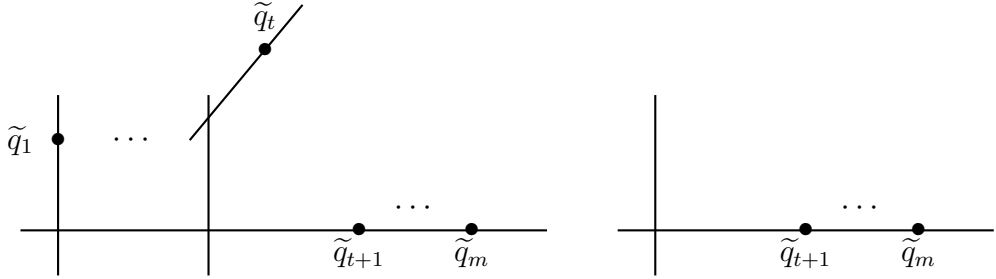
\begin{figure}[h!]
\begin{center}
\begin{tikzpicture}[xscale=0.5,yscale=0.60] [xscale=0.36,yscale=0.36]

       	\draw [thick, black] (-12,0) to (2,0);
	\draw [thick, black] (-11,-1) to (-11,3);
	\node at (-11,2) {$\bullet$};
	\node at (-12,2) {$\widetilde{q}_1$};
	\node at (-9,2) {$\cdots$};
	
\draw [thick, black] (-7,-1) to (-7,3);
\draw [thick, black] (-7.5,2) to (-4.5,5);
\node at (-5.5,4) {$\bullet$};
\node at (-5.5,4.75) {$\widetilde{q}_t$};

\node at (-3,0) {$\bullet$};
\node at (-3,-0.5) {$\widetilde{q}_{t+1}$};
\node at (0,0) {$\bullet$};
\node at (0,-0.5) {$\widetilde{q}_{m}$};
\node at (-1.5,0.5) {$\cdots$};

\end{tikzpicture}
\hspace{0.25in}
\begin{tikzpicture}[xscale=0.5,yscale=0.60] [xscale=0.36,yscale=0.36]

       	\draw [thick, black] (-8,0) to (2,0);
	
\draw [thick, black] (-7,-1) to (-7,3);
		
\node at (-3,0) {$\bullet$};
\node at (-3,-0.5) {$\widetilde{q}_{t+1}$};
\node at (0,0) {$\bullet$};
\node at (0,-0.5) {$\widetilde{q}_{m}$};
\node at (-1.5,0.5) {$\cdots$};

\end{tikzpicture}
\caption{The curves $C'$ (left) and $C''$ (right), obtained by deleting components from $\wt{C}$.}\label{fig:hypersurfacecomb_deletion}
\end{center}
\end{figure}

Because $h^1(C'',f^*T_X|_{C^{\prime\prime}})\le A$, we have that $f|_{C''}:C''\to X$ moves in a family of dimension at most
$$e^{\prime\prime}(n+1-d) + (n-1)(1-g) + A,$$
and $C'$, which is obtained from $C''$ by attaching $t$ free curves of degree $e_i$,
moves in a family of dimension at most
$$ e^{\prime\prime}(n+1-d) + (n-1)(1-g) + A  + \sum_{i=1}^t 
(e_i(n+1-d)-1).$$
Because $C'$ also passes through $m$ general points, we also have
$$ e^{\prime\prime}(n+1-d) + (n-1)(1-g) + A  + \sum_{i=1}^t 
(e_i(n+1-d)-1) \geq m(n-1).$$ 
Since $m(n-1)= e(n+1-d)+(n-1)(1-g),$ and 
$e \geq e^{\prime\prime} +\sum_{i=1}^t e_i$, we get 
$t \leq A$. Therefore for large enough $m$, we have
$m-t \geq g+1$.  We conclude from \cite[Proposition 13 and Proposition 14]{LP23} that $h^1(C,f^*T_X|_C)=0$. In particular, $f|_{C}:C\to X$ moves in a family of the expected dimension.

Now, by \cite{RY16}, the space of rational curves passing through any fixed point of $X$ has the expected dimension, and therefore the space of trees of rational curves passing through any fixed point of $X$ has the expected dimension. As $f|_{C}:C\to X$ also moves in a family of the expected dimension, the same is true of $f$. By assumption, the domain of $f$ is reducible, so $f$ moves in a family of dimension strictly less than $e(n+1-d)+(n-1)(1-g)=m(n-1)$, which therefore cannot dominate $X^{\times m}$. This completes the proof.
\end{proof}

We still need to verify the $h^{1}$ condition used in Lemma \ref{lemm:assumingh1}.  We first need a lemma.

\begin{lemm}\label{chain-rational}
Let $m\ge 3$ be an integer. Let $X\subset\PP^n$ be a smooth hypersurface of degree $d \leq (n+3)/3$, and let
$p_1, \dots, p_m$ in $X$ be general points. Let $(B,q_1,\ldots,q_m)\in \overline{M}_{0,m}$ be the unique point of  $\overline{M}_{0,m}$ if $m=3$ and the pointed stable curve depicted in Figure \ref{fig:B_chain} if $m \geq 4$. So if $m \geq 4$, then $B$ is given by a chain of smooth rational curves $B_1\cup\cdots\cup B_{m-2}$ with $q_i\in B_{i-1}$ for $i=2,\ldots,m-1$ and additionally $q_1\in B_1$ and $q_m\in B_{m-2}$. Assume that there exists an $m$-pointed stable map $f:\wt{B}\to X$ with $f(\wt{q_i})=p_i$ for $i=1,2,\ldots,m$, such that the stable contraction of the domain $(\wt{B},\wt{q_1},\ldots,\wt{q_m})$ is equal to $(B,q_1,\ldots,q_m)$.
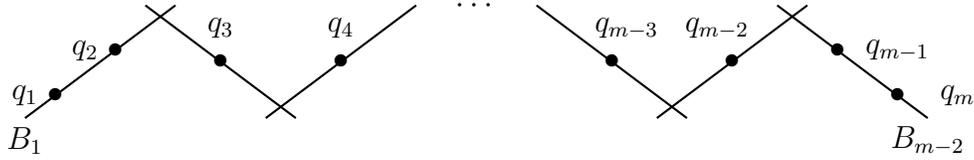
\begin{figure}[h!]
\begin{center}
\begin{tikzpicture}[xscale=0.4,yscale=0.30] [xscale=0.36,yscale=0.36]

       	\draw [thick, black] (-2,0) to (-7,-5);
					\node at (-4.5,-2.5) {$\bullet$};
		\node at (-4.5,-1) {$q_4$};
		\draw [thick, black] (-6,-5) to (-11,0);
		\node at (-8.5,-2.5) {$\bullet$};
		\node at (-8.5,-1) {$q_3$};
		\draw [thick, black] (-10,0) to (-15,-5);
						\node at (-12,-2) {$\bullet$};
				\node at (-13,-2) {$q_{2}$};
				\node at (-14,-4) {$\bullet$};
			\node at (-15,-4) {$q_{1}$};
			\node at (-15,-6) {$B_{1}$};

       \node at (0,0) {$\cdots$};
       
       	\draw [thick, black] (2,0) to (7,-5);
						\node at (4.5,-2.5) {$\bullet$};
		\node at (5,-1) {$q_{m-3}$};
		\draw [thick, black] (6,-5) to (11,0);
				\node at (8.5,-2.5) {$\bullet$};
		\node at (8,-1) {$q_{m-2}$};
		\draw [thick, black] (10,0) to (15,-5);
				\node at (12,-2) {$\bullet$};
				\node at (14,-2) {$q_{m-1}$};
				\node at (14,-4) {$\bullet$};
			\node at (16,-4) {$q_{m}$};
			\node at (15,-6) {$B_{m-2}$};
		
\end{tikzpicture}
\caption{The curve $B=B_1\cup\cdots B_{m-2}$.}\label{fig:B_chain}
\end{center}
\end{figure}

 Then, $e(n+1-d)-(m-1)(n-1) \geq 0$, and if 
$e(n+1-d)-(m-1)(n-1) =0$, then $\wt{B}=B$.
\end{lemm}

\begin{proof}
Let $c:\wt{B}\to B$ be the stabilization map. If $m \geq 4$, then we abusively denote by $B_i$ the unique component of $\wt{B}$ mapping isomorphically to $B_i\subset B$. For $i=1,2,\ldots,m-3$, let $B_{i,i+1}\subset \wt{B}$ be the unique chain of rational curves connecting $B_i$ to $B_{i+1}$ (but not containing either component). Note that $B_{i,i+1}$ may be (and in fact, in the end, will be) empty. For $i=1,2,\ldots,m-3$, let $S_i=B_i\cup B_{i,i+1}$, and let $S_{m-2}=B_{m-2}$. Let $\wt{B}_s\subset\wt{B}$ be the union of all of the $S_i$. A piece of the curve $\wt{B}$ is depicted in Figure \ref{fig:Btilde_chain}.
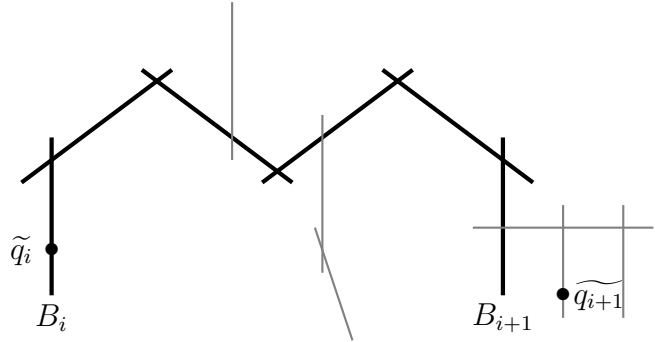
\begin{figure}
\begin{center}
\begin{tikzpicture}[xscale=0.4,yscale=0.30] [xscale=0.36,yscale=0.36]

       	\draw [ultra thick, black] (-14,-3) to (-14,-10);
	\node at (-15,-8) {$\wt{q_i}$};
	\node at (-14,-8) {$\bullet$};
	\node at (-14,-11) {$B_i$};
	
	\draw [ultra thick, black] (-10,0) to (-15,-5);
	
	\draw [ultra thick, black] (-6,-5) to (-11,0);
	\draw [thick, gray] (-8,-4) to (-8,3);
	
       	\draw [ultra thick, black] (-2,0) to (-7,-5);
		\draw [thick, gray] (-5,-2) to (-5,-9);
		\draw [thick, gray] (-5.25,-7) to (-4,-12);
	
       	\draw [ultra thick, black] (-3,0) to (2,-5);
	
	       	\draw [ultra thick, black] (1,-3) to (1,-10);
		\node at (1,-11) {$B_{i+1}$};
		\draw [thick, gray] (0,-7) to (6,-7);
		\draw[thick, gray] (3,-6) to (3,-11);
		\node at (3,-10) {$\bullet$};
		\node at (4.2,-10) {$\wt{q_{i+1}}$};
		\draw[thick, gray] (5,-6) to (5,-11);

\end{tikzpicture}
\caption{A piece of the curve $\wt{B}$. Here, the subcurve $B_{i,i+1}$ is given by the union of the four components connecting $B_i$ to $B_{i+1}$, and $S_i$ is the chain of five curves $B_i\cup B_{i,i+1}$. The bold components are part of the spine $\wt{B_s}$, whereas the gray components are among the trees $T$. }
\label{fig:Btilde_chain}
\end{center}
\end{figure}

Note that $\wt{B}_s\subset\wt{B}$ is itself a chain of rational curves, and $\wt{B}$ is obtained from $\wt{B}_s$ by attaching pairwise disjoint trees of rational curves at smooth points. Let $T$ be such a rational tree. The stability of $f$ shows that $T$ cannot be contracted by $f$. 

Let $\wt{B}_i$ be the union of $S_i$ and all of the rational trees $T$ attached to $S_i$. Thus, $\wt{B}$ is the union of all $\wt{B}_i$, for $i=1,2,\ldots,m-2$, and every component of $\wt{B}$ belongs to exactly one $\wt{B}_i$.

We now proceed to the proof of the Lemma.  We may assume that $e(n+1-d) \leq (m-1)(n-1)$; we argue by induction on $m$ that equality must hold (and the resulting curve has the claimed form). If $m=3$ then let $\wt{B}_s$ be the unique component of $\wt{B}$ mapping isomorphically to $B$. 
Then $e\le\frac{2(n-1)}{n+1-d}$ and since $d \leq (n+3)/3$ we have $e\le 2$.
So we have the following cases: 
\begin{enumerate}
\item $\deg(\wt{B}_s)=2$, in which case $\wt{B}$ must be irreducible,
\item $\deg(\wt{B}_s)=1$, in which case $\wt{B}_s$ must map to a line between two of $p_1,p_2,p_3$, and a second component $T\subset\wt{B}$ must map to a line through the third point, or
\item $\deg(\wt{B}_s)=0$, in which case $\wt{B}_s$ must be contracted to one of $p_1,p_2,p_3$, 
and two other components must be attached to $\wt{B}_s$ at the other two points, mapping isomorphically to lines on $X$.
\end{enumerate}
In the first case, $f:\wt{B}\to X$ is an irreducible free conic, thus moving in a family of the expected dimension, and such an $f$ passes through 3 general points only when $e(n+1-d)=(m-1)(n-1)$.   (Here, we get simply $d=2$, that is, $X$ is a quadric.) In the second case, $f$ is a chain of two free lines, which moves (as a pointed stable map) in a family of the expected dimension of $e(n+1-d)+(n-1)-1< m(n-1)$,
so $f$ cannot pass through $m=3$ general points. The third case is similarly impossible. This establishes the base case.

Suppose now that the Lemma holds for every $m'<m$.  We use the induction hypotheses and the assumption $d\leq (n+3)/3$ to show that $\wt{B}$ is a chain of lines and free conics to which some free lines are attached. Fix an index $i$ and let $\wt{B}'=\bigcup_{j<i}\wt{B_j}$ and
$\wt{B}^{\prime\prime}=\bigcup_{j>i}\wt{B_j}$. Note that $\wt{B}'$ or $\wt{B}^{\prime\prime}$ is empty if $i=1$ or $i=m-2$, respectively. Let $e',e^{\prime\prime}$ be the total degree of $f$ restricted to $\wt{B}',\wt{B}^{\prime\prime}$, respectively, and let $e_i$ be the degree upon restriction to $\wt{B_i}$.

Let $a' =  e'(n+1-d)-(i-1)(n-1)$ 
and $a^{\prime\prime} =  e^{\prime\prime}(n+1-d)-((m-i-1)-1)(n-1)$. If $3 \leq i \leq m-2$, then by our induction hypothesis $a' \geq 0$. If $i=2$, then either $d=1$, or 
$d \geq 2$ and there is no line through two general points of $X$, so $e' \geq 2$. So in all cases $a' \geq 0$. Similarly $a^{\prime\prime} \geq 0$. Because
$e=e'+e^{\prime\prime}+e_i$, we have
\begin{align*}
a'+a^{\prime\prime} &=e(n+1-d)-(m-1)(n-1) -e_i(n+1-d)+2(n-1) \\
&\leq -e_i(n+1-d)+2(n-1).
\end{align*}
If $e_i(n+1-d)>2(n-1),$ we obtain a contradiction. Therefore, we have
$e_i \leq \frac{2(n-1)}{n+1-d} <3$, so $e_i\leq 2$. 

We now reduce to the case in which every component of $\wt{B}$ that is a leaf (that is, is connected to only one other component of $\wt{B}$) contains a marked point. Indeed, if $\wt{B}$ contains a leaf without a marked point, then deleting that component, and stabilizing $f$ if necessary, decreases $e$ while leaving $m$ the same. Thus, it suffices to show that the new stable map, of strictly smaller degree, cannot exist. 

In particular, any tree attached to $\wt{B}_s$ must contain a marked point. By construction, this marked point is unique, so in fact, the tree must be a chain, and all components must map to $X$ with positive degree.

Now, let $s_i=\deg(S_i)$ and $t_i=\deg(\wt{B}_i)-\deg(S_i)$.  Since $e_{i} \leq 2$, the only possibilities are that:
\begin{enumerate}
\item $t_i=0$, in which case $s_i\le 2$,
\item $t_i=1$, in which case $s_i\le 1$, or
\item $t_i=2$, in which case $s_i=0$.
\end{enumerate}
Case (3) is only possible if $i=1$ or $i=m-2$. Indeed, if $2\le i\le m-3$, then after contracting $S_i$ and deleting the (necessarily unique) tree attached to it, we obtain an instance of this Lemma with $e$ replaced by $e-2$ and $m$ replaced by $m-1$. In this case, since $d \leq (n+3)/3$, we have $(e-2)(n+1-d)-(m-2)(n-1) <0$, so $f$ cannot exist by the inductive hypothesis. On the other hand, if $i=1$ or $i=m-2$, then the same argument and the stability of $f$ force $S_i$ to consist of a single contracted component, with two free lines passing through the two points $p_1,p_2$ or $p_{m-1},p_{m}$.

In case (2), the unique tree attached to $S_i$ must consist of a single line containing a marked point, so in particular, the line is free. In case (1), $S_i$ consists of a 
free irreducible conic or a union of two lines. 

To conclude, $\wt{B}$ consists of a chain of components of degree 0, 1, or free components of degree 2 (and the contracted components are necessarily stable), with possibly some additional free lines containing marked points $\wt{q_i}$ attached to the chain. Note that after contracting the degree $0$ components in $\wt{B}_s$, we get a chain of lines and free conics, so by 
Proposition \ref{chain_conics}, $\wt{B}_s$ moves in a family of expected dimension. Since the family of free lines through any point of $X$ has the expected dimension $n-d-1$, it follows that $f$ moves in a family of the expected dimension, and can therefore only pass through all of the $p_i$ if $e(n+1-d)-(m-1)(n-1) =0$ and $\wt{B}=B$.
\end{proof}

We first prove the genus $0$ case of Theorem \ref{hypersurface-main} using the previous lemmas.  We then build off this case to prove the statement for arbitrary genus.

\begin{prop}\label{prop:h1_bounded}
Let $X$ be a smooth hypersurface of degree $d$ in $\PP^n$ as in Theorem \ref{hypersurface-main}.
 For any integer $a \geq 0$, there is a constant 
$Q_a$ with the following property:
if $e, m$ are such that  $e(n+1-d)-(m-1)(n-1) \leq a$, then for a general $(C, p_1, \dots, p_m) \in \overline{M}_{0,m} \times X^{\times m}$ and for every $(\wt{C}, f, \widetilde{q}_1,\dots, \widetilde{q}_m)$ in the fiber of $\phi$ over $(C,p_1, \dots, p_m)$, we have $h^1(\wt{C},f^*T_X) \leq Q_a$. 
\end{prop}

\begin{proof}

We prove the statement by induction on $a$. First suppose $a=0$. If $X$ is general, then the statement follows from the irreducibility of the space of rational curves of any given degree on a general hypersurface of degree $\leq n-2$ in $\PP^n$. 

Now assume $a=0$ and $X$ is an arbitrary  smooth hypersurface of degree 
$d  \leq (n+3)/3$ in $\PP^n$. Then by Proposition~\ref{chain_conics} the spaces of chains of lines and free conics in $X$ have the expected dimension. 
We show for general $[C] \in \overline{M}_{0,m}$ and for general points $p_1, \dots, p_m \in X$, we have $h^1(\wt{C},f^*T_X)=0$ for every $[f]$ in the fiber of $\phi$ over $([C],p_1, \dots, p_m)$, so $Q_0=0$. By upper semicontinuity of $h^1$, it is enough to show there exists $(C,p_1, \dots, p_m) \in \overline{M}_{0,m} \times X^{\times m}$ such that $h^1(\wt{C},f^*T_X)=0$ for every point $[f]$ in the fiber of $\phi$ over $(C,p_1, \dots, p_m)$.  Taking $C$ to be the point $(B,q_1,\ldots,q_m)\in \overline{M}_{0,m}$ defined in the statement of Lemma \ref{chain-rational}, the conclusion of the same Lemma shows that $C$ has the required property.

   Next suppose the statement holds for any integer smaller than $a$. Then for general $C$, $p_1, \dots, p_m$, we use the same argument as in Lemma \ref{chain-rational}. We define $C_i,C_{i,i+1},S_i,\wt{C}_s,\wt{C}_i\subset\wt{C}$ as in the proof of Lemma \ref{chain-rational} (where we have replaced all instances of $B$ with $C$), and for any index $i$, define $e', e^{\prime\prime},a', a^{\prime\prime},e_i$, as before.  Then, 
    $$a'+a^{\prime\prime} = a -e_i(n+1-d)+2(n-1).$$
    
Thus, if $e_i \geq 3$ for some $i$, then since $a',a^{\prime\prime} \geq 0$, we have $a', a^{\prime\prime} <a$. Also, $e_i$ is bounded in terms of $a, n,d$ (since for example $a-e_i(n+1-d)+2(n-1) \geq 0$).  Thus there is a constant $A$ which depends only on $a, n,d$ such that the contribution to
    $h^1$ from $\wt{C_i}$ is at most $A$. Therefore $h^1(\wt{C},f^*T_X)\le A+Q_{a'}+Q_{a''}+2n-2,$ and we are done. 
    
    So we may assume $e_i \leq 2$ for every $i$. As in the proof of Lemma \ref{chain-rational}, we next reduce to the case in which all leaves of $\wt{C}$ contain a marked point. A leaf without a marked point would need to have degree 1 or 2, and deleting it from $\wt{C}$ and stabilizing the resulting stable map gives an instance of the Proposition with $a$ replaced by some $a'<a$. Then, $h^1(\wt{C},f^{*}T_X)\le Q_{a'}+A+(n-1)$,
 where $A$ denotes here the maximum value of $h^1(f^{*}T_X)$ for a line or conic on $X$. As this operation can be performed no more than $a$ times, we obtain an upper bound on $h^1(\wt{C},f^{*}T_X)$ depending only on $a$.
    
Now, we have the same three cases for the degrees of $s_i,t_i$ as in Lemma \ref{chain-rational}. As before, if we are in case (3) and there exists a tree of degree at least 2 attached to $\wt{C}_s$, then deleting it and contracting $S_i$ yields an  instance of the Proposition with $a$ replaced by some $a'<a$. As before, $h^1(\wt{C},f^{*}T_X)$ is bounded above in terms of (the finitely many constants) $Q_{a'}$, $n$, and the maximum value of $h^1(f^{*}T_X)$ for a line or conic on $X$.
     
Finally, we conclude that otherwise, as in Lemma \ref{chain-rational}, $f$ is given by a chain of lines and free conics, possibly along with free lines attached as tails, and contracted components. In particular, $f$ moves in a family of the expected dimension.  Note that all conics comprising $f$ are free, as they arise from case (1) of the analysis of Lemma \ref{chain-rational}. 

On the other hand, $f$ may contain non-free lines $L$, which necessarily do not contain any of the marked points. We claim, however, that there can be at most $a$ such components. Indeed, each additional component with no marked point decreases the expected dimension, hence the actual dimension by $1$. By assumption, the expected dimension of the space of maps is $a$ more than the dimension of $\cM_{g,m}\times X^{\times m}$. So if there are more than $a$ non-free components, the space of maps cannot dominate $\cM_{g,m}\times X^{\times m}$.

Now, the only non-trivial contributions to $h^1(\wt{C},f^{*}T_X)$ come from such components $L$, of which there are a bounded number.  Note that $h^1(L,T_X|_L)$ is uniformly bounded depending only on $X$. It follows again that $h^1(\wt{C},f^{*}T_X)$ is bounded uniformly in terms of $X$ and $a$. 
\end{proof}

\begin{proof}[Proof of Theorem \ref{hypersurface-main} for arbitrary $g$] 
It suffices to verify the hypothesis of Lemma \ref{lemm:assumingh1}. Assume that no such constant $A$ as in the statement of the Lemma exists. Then, given any integer $A'$, there exist integers $m$ and $e$ satisfying $e(n+1-d)=(m+g-1)(n-1)$ and with the following property. Let $[C]=(C, q_1, \dots, q_m)$ be a point of $\overline{\mathcal{M}}_{g,m}$ given by a union of a general pointed curve $(C_0,q_1,\ldots,q_m)\in M_{0,m}$ and a smooth curve $C_g$ of genus $g$, attached at a smooth point of $C_0$, see Figure \ref{fig:genus_g_to_genus_0}. Then, there exists $(\widetilde{C}, f, \widetilde{q}_1, \dots, \widetilde{q}_m)$ in the fiber over $([C],p_1,\ldots,p_m)\in \mathcal{M}_{g,m}\times X^{\times m}$ such that $h^1(\widetilde{C},f^*T_X)>A'$.

\begin{figure}[h!]
\begin{center}
\begin{tikzpicture}[xscale=0.4,yscale=0.30] [xscale=0.36,yscale=0.36]
			
        \node at (-13,-6) {$C_g$};
         \node at (13,-6) {$C_0$};
	\draw [thick, black] (-12,-6) to (2,1);
	\draw [thick, black] (12,-6) to (-2,1);
	\node at (4,-0.5) {$q_1$};
		\node at (4,-2) {$\bullet$};
			\node at (6,-1.5) {$q_2$};
		\node at (6,-3) {$\bullet$};
				\node at (8,-3) {$\cdots$};
				\node at (10,-3.5) {$q_n$};
		\node at (10,-5) {$\bullet$};
		
\end{tikzpicture}
\caption{The curve $C=C_g\cup C_0$ with $q_i\in C_0$.}\label{fig:genus_g_to_genus_0}
\end{center}
\end{figure}
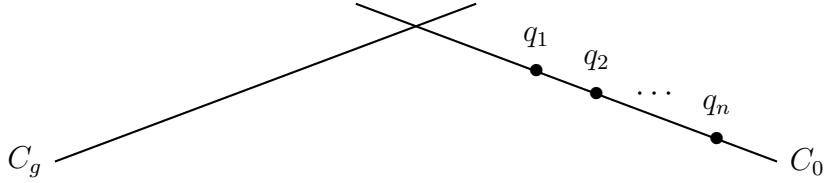

Write $\widetilde{C}=\widetilde{C}'\cup \widetilde{C}^{\prime\prime}$, so that $\widetilde{C}'$ is a tree of rational curves
and contains all $m$ marked points, and $\widetilde{C}^{\prime\prime}$ contains $C_g$, possibly along with additional rational components.  Let $e',e''$ be the total degrees of the 
restriction of $f$ to $\widetilde{C}', \widetilde{C}^{\prime\prime}$, respectively.

Then
$$e'(n+1-d)-(m-1)(n-1)\leq e(n+1-d)-
(m-1)(n-1)=g(n-1),$$ so by Proposition \ref{prop:h1_bounded}, $h^1(f^*T_X|_{\widetilde{C}'})$ is bounded by a constant. On the other hand, if $$e''>\frac{g(n-1)}{n+1-d}=e-\frac{(m-1)(n-1)}{n+1-d},$$ then
$$e'(n+1-d)-(m-1)(n-1)=(e-e'')(n+1-d)-(m-1)(n-1)<0,$$
so the case $g=0$ implies that no such $f$ can exist. Therefore, $e''$ is bounded, so $h^1(f^*T_X|_{\widetilde{C}''})$ is bounded as well. 

On the other hand, we have assumed that $h^1(\widetilde{C},f^*T_X)>A'$ for some $A'$ arbitrarily large, so we have reached a contradiction.  The conclusion now follows from Lemma \ref{lemm:assumingh1}.
\end{proof}

\begin{rema}
The same induction argument as in the  proofs of  Lemma  \ref{chain-rational} and Proposition \ref{prop:h1_bounded} shows that when $d<(n+3)/2$, we may only consider the case where $e_i \leq 3$ for every $i$. So the statements, and hence the statement of Theorem \ref{theo:introhypersurface},  hold for any smooth hypersurface $X$ in this degree range  if we know families of certain trees of lines, conics, and free cubics on $X$ have the expected dimension. In fact, one can use Bend and Break to reduce this statement to a question about families of trees of lines.

\end{rema}

\bibliographystyle{alpha}
\bibliography{asymenum}

\end{document}